\documentclass[preprint,12pt]{elsarticle}

\journal{Computer-Aided Design}

\usepackage[margin = 2.5cm]{geometry}
\usepackage[T1]{fontenc}
\usepackage[utf8]{inputenc}
\usepackage[table]{xcolor}
\usepackage{times,amsmath,amsfonts,amssymb,mathtools,marvosym,algorithm, algorithmicx,algpseudocode,multirow,etoolbox}
\usepackage[binary-units=true]{siunitx}
\usepackage[inline]{enumitem}
\usepackage{tikz}
\usepackage[english]{babel}
\usepackage{hyperref}
\usetikzlibrary{calc}

\newcommand{\TITEL}{Multilevel T-spline Approximation for Scattered Observations with Application to Land Remote Sensing}

\renewenvironment{abstract} {\begin{quotation} \textbf{Abstract.}\quad} {\end{quotation}}

\hypersetup{
    pdftitle = {\TITEL},
    pdfauthor = {Ga\"el Kermarrec, Philipp Morgenstern},
    pdfkeywords = {T-splines} {Multilevel B-spline Approximation} {NURBS} {surface approximation} 
    {adaptive local refinement},
    colorlinks=true,                 
    linkcolor=red!50!black,          
    citecolor=green!50!black,        
    filecolor=magenta!50!black,      
    urlcolor=cyan!50!black           
}

\definecolor{lightgray}{gray}{0.925}
\definecolor{headrow}{gray}{0.7}
\let\oldtabular\tabular
\let\endoldtabular\endtabular
\renewenvironment{tabular}{\rowcolors{2}{lightgray}{white}\oldtabular}{\endoldtabular}

\newlist{inlinelist}{enumerate*}{1}
\setlist[inlinelist]{label=\emph{(\roman*)}, itemjoin={{, }}, itemjoin*={{, and }}}

\newcommand{\R}{\mathbb R}
\newcommand{\mesh}{\mathcal Q}
\newcommand{\nodes}{\mathcal N}
\newcommand{\knots}{\mathsf k}
\newcommand{\degree}{p}
\newcommand{\nobs}{{n_\text{obs}}}
\newcommand{\thr}{\mathsf{T\!H}}
\newcommand{\rmse}{\text{RMSE}}
\newcommand{\maxerr}[1][]{\operatorname{Max\_err}_{#1}}
\newcommand{\ct}{\mathsf{C\!T}}

\newcommand{\Ssmooth}{S_\text{smooth}}
\newcommand{\Soutliers}{S_\text{smooth}^\text{outliers}}
\newcommand{\Sgaps}{S_\text{smooth}^\text{gap}}
\newcommand{\Ssharp}{S_\text{sharp}}
\newcommand{\tSsmooth}{\tilde S_\text{smooth}}
\newcommand{\tSoutliers}{\tilde S_\text{smooth}^\text{outliers}}
\newcommand{\tSgaps}{\tilde S_\text{smooth}^\text{gap}}
\newcommand{\tSsharp}{\tilde S_\text{sharp}}
\newcommand{\ncpk}{n_{\text{cp},k}}
\newcommand{\rmsemathk}[1][k]{\rmse_{\text{math},#1}}
\newcommand{\rmsenoisek}[1][k]{\rmse_{\text{noise},#1}}
\DeclareMathOperator{\supp}{supp}
\makeatletter

\newcommand*{\transp}{%
  ^{\mathpalette\@transpose{}}%
}
\newcommand*{\@transpose}[2]{%
  \raisebox{\depth}{$\m@th#1\intercal$}%
}
\makeatother

\begin{document}
\begin{frontmatter}

\title{\TITEL}

\author[gis]{Gaël Kermarrec\corref{cor}}
\ead{kermarrec@gih.uni-hannover.de}%

\author[ifam]{Philipp Morgenstern}
\ead{morgenstern@ifam.uni-hannover.de}

\address[gis]{Geodetic Institute, Leibniz Universität Hannover, Nienburger Str. 1, 30167 Hannover, Germany}
\address[ifam]{Institute of Applied Mathematics, Leibniz Universität Hannover, Welfengarten 1, 30167 Hannover, Germany}
\cortext[cor]{\rule{0pt}{2em}corresponding author}

\begin{abstract}
 In this contribution, we introduce a multilevel approximation method with T-splines for fitting scattered point clouds iteratively, with an application to land remote sensing.
This new procedure provides a local surface approximation by an explicit computation of the control points and is called a multilevel T-splines approximation (MTA). It is computationally efficient compared with the traditional global least-squares (LS) approach, which may fail when there is an unfavourable point density from a given refinement level. 
We validate our approach within a simulated framework and apply it to two real datasets:
\begin{inlinelist}
 \item a surface with holes scanned with a terrestrial laser scanner  
 \item a patch on a sand-dune in the Netherlands
\end{inlinelist}.
Both examples highlight the potential of the MTA for rapidly fitting large and noisy point clouds with variable point density and with similar results compared to the global LS approximation.
\end{abstract}



\begin{keyword}
T-splines \sep multilevel B-splines approximation \sep  NURBS \sep  surface approximation \sep  adaptive local refinement \sep GIS
\PACS 02.60.Ed
\MSC 33F05
\end{keyword}

\end{frontmatter}

\section{Introduction}\label{sec: intro}

Parametric spline surface approximation is an important step in reverse engineering to convert real object data to a computer-aided design (CAD) model \cite{RajaFernandes2008}. The observations can be represented in a more structured and condensed format, which is favourable for further processing. Similar approximation techniques can be applied to point clouds representing land, vegetation, architectural or engineering structures (dams, bridges or tunnels). Deformation analysis based on rigorous statistical testing is a prominent application \cite{KermarrecKargollAlkhatib2020,HarmeningNeuner2020}.
The initial stage is usually the scanning of an object or scene. A contactless sensor, for example, a terrestrial laser scanner (TLS), can be used for this.
This scanner is widely employed for high-resolution terrestrial monitoring within a geographic information system (GIS) context: large and scattered point clouds are obtained which may be unstructured and with variable data density. This particularity challenges the usual spline-based surface approximation techniques, such as those used for CAD. It makes specific improvements of the fitting method mandatory. The latter can be divided into 
\begin{inlinelist}
 \item non-adaptive methods, for which the optimal surface is globally adapted in the approximation process
 \item adaptive methods
\end{inlinelist}.
The fitting strategy commonly used utilising B-splines belongs to the class of non-adaptive methods, as presented in \cite{FarrRosenCaroAndManyMore2007}, with observations from the Shuttle Radar Topography Mission, or \cite{BureickAlkhatibNeumann2016,Koch2009} for TLS observations.
This approximation method is greatly limited by the tensor product formulation of non uniform rational basis splines (NURBS). The degrees of freedom of the approximation become unnecessarily large in regions where refinement was not indicated, which can lead potentially to approximating the noise, i.e. \emph{overfitting}.
It is advantageous to consider approaches with local mesh refinement for geodetic applications, such as deformation analysis between point clouds \cite{ShapiroBlaschko2004}, which are often large, scattered and noisy (e.g.\ for TLS, range noise \cite{WujanzBurgerMettenleiterNeitzel2017,KermarrecLoeslerHartmann2021} and angle noise \cite{KermarrecHartmann2021}). They avoid the risk of overfitting to some extent and are more efficient in terms of computation time (CT) and the number of parameters to estimate. 

There are several alternative approaches for locally refined and linearly independent B-splines in two and three dimensions. We mention here
hierarchical B-splines (HB-splines \cite{ForseyBartels1988, Kraft1997}) and their truncated version (THB-splines \cite{GiannelliJuettlerSpeleers2012}),
locally refined B-splines (LR B-splines \cite{DokkenLychePettersen2013,PatriziManniPelosiSpeleers2020}),
polynomial splines over T-meshes (PHT-splines \cite{LiZhang2010,SchumakerWang2012}),
and analysis-suitable T-splines (AST-splines \cite{SederbergZhengBakenovNasri2003,LiZhengSederbergHughesScott2012,ScottLiSederbergHughes2012}).
All these approaches have been applied successfully to surface fitting \cite{Koch2009,ForseyBartels1996,BraccoGiannelliGrossmannSestini2018,SkyttHarphamDokkenDahl2017,NiWangDeng2019,Wang2009}. In this contribution,
we propose to use the bicubic AST-splines.
A  comparison with other refinement strategies is beyond the scope of this paper and is left to further contributions. Interested readers can refer to the study of \cite{Remonato2014,HennigKaestnerMorgensternPeterseim2017}.

Concerning the fitting of a scattered point cloud, the $z$ component of the point cloud is parametrized regarding the $(x,y)$ co-ordinates (see \cite{FloaterHormann2005}).
The parametrized point cloud is then approximated by a T-spline surface. Starting from a coarse mesh, successive refinements are performed in domains where the distance between the data points and the approximated surface exceeds a given threshold. 
The approximation of the point cloud is performed by least-squares (LS) adjustment: the $L^2$ norm between observations and the parametric surface is minimized. 
This approach is called surface skinning for gridded data \cite{Woodward1988} and was used to reconstruct a spline surface from $z$-map data \cite{ZhengWangSeah2005}.  
The aim of the LS approximation is to determine a set of coefficients or control points that yields a sufficiently close approximation of the given point cloud by spline functions regarding the $L^2$ norm. A fairness functional is often included to smooth and stabilize the approximation \cite{CelnikerGossard2017}.
The system of equations can be solved using a preconditioned conjugate gradient method, as in \cite{SkyttBarrowcloughDokken2015} or a QR solver 
\cite{Chen1989}%
. Unfortunately, because of the variable data density in real point clouds, many strategies may reach their limits, even for linearly independent splines and in the absence of data gaps \cite{BraccoGiannelliGrossmannSestini2018}. 
As the refinement level increases, iterative approaches for solving the LS problem may become very slow and inaccurate. This challenge is enhanced by the addition of a smoothing term; \cite{DokkenLychePettersen2013} showed exemplarily that the condition number of the stiffness matrix for HB, truncated HB and LR B-splines is linked to increasing refinement levels in the context of isogeometric analysis.
An adaptive procedure was proposed recently by \cite{WangZouAndOthers2021} based on a local LS fitting method to face that challenge. Unfortunately, the fitting of complex geometries still produced ripples in the presence of noise and outliers. With the aim of speeding up the algorithm and improving the fitting in terms of quality and conciseness, \cite{FengTaguchi2017} developed a novel split-connect-fit algorithm. 
This method is based on the division of the point cloud into a set of B-spline patches that are connected into a single T-spline surface.
It may not be applicable for point clouds with non-uniform parametrization.
Furthermore, the proposal is still based on a refinement using a conjugate gradient method, which is expected to remain slow as the refinement becomes more sophisticated. 
In the same direction, \cite{Lu2020AFT} introduced a segmentation technique to identify the inactive and active regions of the T-grid, and proposed a fast T-spline fitting method relying on T-grid segmentation. 
\cite{FuEtAl2015} proposed a T-spline surface reconstruction method based on the centroidal Voronoi tessellation sampling strategy using constrained optimization theory. 
The Gaussian curvature of the surface was taken as the regional density function in the parameter region of the originally designed NURBS surface, and the centroidal Voronoi tessellation was generated by the Lloyd algorithm. 
The fitting methods for all those proposals often 
\begin{inlinelist}[itemjoin*={{, or }}]
\item still require a large number of control points that need to be estimated
\item are based on patches that need to be fused together using rather challenging continuity options
\item rely on complex mathematical developments, making further implementations and applications in usual CAD or GIS software not straightforward
\end{inlinelist}. 
Therefore, a comparison with these approaches is beyond the scope of this paper.
Additionally, the major problem of variable data density is only partially solved: the need for an easy-to-use, computationally low-cost, local and data-dependant approximation based on an effective mesh refinement strategy is still present.  

Complementary to the LS method, the multilevel B-spline approximation (MBA) was introduced by \cite{LeeWolbergShin1997} as a local approximation method. The MBA is an explicit procedure for which no equation system needs to be solved. It has been used for scalar data approximation in geology and oceanography \cite{HaberZeilfelderDavydovSeidel2008} and for three-dimensional (3D) object reconstruction \cite{LeeWolbergShin1997}. The MBA has some similarities to quasi-interpolants (see, e.g. \cite{SpeleersManni2016}, for an application to hierarchical spline spaces). However, the MBA coefficients are not independent of each other as there is a certain ordering according to the levels. The MBA cannot, thus, be rigorously classed as a quasi-interpolant. The algorithm for the MBA starts by fitting a B-spline surface with a coarse mesh to the observations. The residuals of the data points obtained from the last fitted surface are recursively approximated using finer meshes, which saves storage. This method was empirically shown to be particularly suitable when the observations are scattered, noisy and non-gridded and/or when data gaps or outliers occur: it is more likely with the MBA that the “right” level of resolution and, therefore, the “right” number of degrees of freedom is chosen, consequently avoiding ripples in the approximated surface. 
The MBA was extended to the HB and LR frameworks by \cite{ZhangTangLi1998} and \cite{SkyttBarrowcloughDokken2015}, respectively, but has not yet been derived for T-splines. This contribution aims to fill this gap. 
More specifically, we will show how
\begin{itemize}
 \item the T-mesh properties of the refinement strategy proposed can be exploited within the MBA. 
The T-mesh is a key component of the adaptive approximation because its refinement enlarges the approximation space in each iteration step. Different methods for refinement exist; in this contribution, we chose the one developed by \cite{MorgensternPeterseim2015}, which ensures an analysis-suitable mesh for which the linear independence of T-splines is given. 
 \item Simulated  point clouds with data gaps and outliers can be approximated efficiently by combining LS, the MBA and the T-splines adaptive refinement mentioned previously, i.e. validating the new methodology.
 \item Observations from land remote sensing from a sand-dune scanned in the Netherlands can be optimally fitted for deformation analysis. This application aims to democratize the use of T-splines for new applications within a GIS context, i.e. not only restricted to CAD.
\end{itemize}
We call our method “MTA” and the traditional LS procedure “LS-T”. 
The remainder of this paper is as follows: in \autoref{sec: maths}, we briefly explain the construction of
NURBS, T-splines, the corresponding geometry map, LS approximation and the MTA.
In \autoref{sec: simulations}, we validate the MTA as a complementary alternative to LS-T with several simulated point clouds. Furthermore, in \autoref{sec: real data}, we show the potential of the MTA for approximating real point clouds when the data density challenges the LS-T.


\section{Mathematical derivation}\label{sec: maths}
In this section, we address the drawbacks of the tensor product structure of NURBS by introducing T-splines. We will further develop the principles of adaptive surface fitting with LS-T and introduce the MTA method as an explicit procedure to approximate a surface with T-splines.

\subsection{Geometry map}
We describe a mapping $S:D\to\R^3$ from a rectangular domain $D\subset\R^2$ to a curved surface in $\R^3$. We call the rectangular preimage $D$, the \emph{parametric domain} and $\R^2$ the \emph{parametric space} with parametric direction $u$ and $v$; $S(D)$ and $\R^3$ the \emph{physical domain} and \emph{physical space}, respectively, with physical directions $x$, $y$ and $z$. 
We define a set of spline functions $B_i(u,v)$, $j=1,\dots,n$ below, which forms a convex partition of unity.
Together with a set of control points $P_j$, $j=1,\dots,n$, this defines a surface via $S(u,v)=\sum_{j=1}^nB_j(u,v)P_j$, hence each point in $D$ is mapped to a linear combination of the control points $P_1,\dots,P_n$. The smoothness of the functions $B_j$ yields a smooth surface as the image of $S$, and their positivity yields that the modelled surface is in the convex hull of the control points.

\subsection{NURBS surface}
In the following, we call a knot vector a non-decreasing sequence of coordinates in the parameter space, where each coordinate is called a knot. A knot span is the interval bounded by two adjacent coordinates. 

Given a polynomial degree $\degree\in\mathbb N$, a 
number of basis functions $n$, and 
a knot vector $\knots=\{u_1,\dots,u_{n+\degree+1}\}$,
a set of univariate B-splines $N_{1,\degree}^\knots,\dots,N_{n,\degree}^\knots$ 
of degree $\degree$ is defined by
the following recursion (known as the Cox-de-Boor formula, \cite{Boor1971}):
\begin{alignat}{2}
 N_{i,0}^\knots(u) &= \begin{cases} 1 & u_i\le u<u_{i+1} \\ 0 & \text{otherwise} \end{cases}
&\quad&\text{for }i=1,\dots,n+\degree \\
 N_{i,q}^\knots(u) &= \frac{u-u_i}{u_{i+q}-u_i}N_{i,q-1}^\knots(u)+\frac{u_{i+q+1}-u}{u_{i+q+1}-u_{i+1}}N_{i+1,q-1}^\knots(u)
&\quad&\begin{multlined}\text{for }q=1,\dots,\degree \\ \text{and }i=1,\dots,n+\degree-q\end{multlined}%
\label{eq: cox de boor}%
\end{alignat}%
These basis functions have important properties, such as the partition of unity, positivity, compact support and continuity. 
In the case of a bicubic NURBS surface with $\degree=3$, the knot information is given by two knot vectors $\knots_u$ and $\knots_v$ of length $n+4$ and $m+4$, respectively, usually ranging from 0 to 1. 
The reader is referred to \cite{FloaterHormann2005} for a review of parametrization techniques.

Let $B_i^u=N_{i,3}^{\knots_u}$ be the $i$-th cubic B-spline from \autoref{eq: cox de boor} above applied with the knot vector $\knots_u$, and $B_j^v=N_{j,3}^{\knots_v}$  the $j$-th cubic B-spline from \autoref{eq: cox de boor} applied with the knot vector $\knots_v$. Given control points $P_{ij}$ and associated weights $w_{ij}$, the NURBS surface is defined by a tensor product of B-spline functions as 
\begin{align}
S(u,v)&=\sum_{i=1}^n\sum_{j=1}^mB_{ij}(u,v)P_{ij}
\\\text{with}\quad
B_{ij}(u,v) &=\frac{w_{ij}B_i^u(u)B_j^v(v)}{\sum_{k=1}^n\sum_{\ell=1}^mw_{k\ell}B_k^u(u)B_\ell^v(v)}
 \quad\begin{multlined}
 \text{for }i=1,\dots,n \\\text{ and }j=1,\dots,m
 \end{multlined}
\end{align}

\subsection{T-splines}\label{sec: tsplines}
\subsubsection{Definition}
T-splines are one realization of B-splines on locally refined meshes, beneath other approaches such as HB-splines, THB-splines and LR-splines outlined in \autoref{sec: intro}. We explain below the construction of bicubic T-splines. For other polynomial degrees, see \cite{VeigaBuffaSangalliVazquez2013}, or \cite{GoermerMorgenstern2021a} for multivariate T-splines.
We suppose that we are given a mesh $\mesh$ that consists of rectangles which are axis-parallel, have pairwise disjoint interiors and their union is the parametric domain. This mesh is called a T-mesh, and we emphasize that it may contain T-junctions. 
We refer to the set of nodes or vertices of the mesh (i.e.\ the set of all corner points of all elements of $\mesh$) as $\nodes(\mesh)$, or $\nodes$ if there is only one mesh considered in the context.
We associate to each vertex $z\in\nodes$ a knot vector in each axis direction, $\knots_u(z)=(u_{-2},\dots,u_2)$ and $\knots_v(z)=(v_{-2},\dots,v_2)$, where the coordinates of $z=(u_0,v_0)$ are the middle entries and the other entries are the $u$-coordinates of the nearest vertical and the $v$-coordinates of the nearest horizontal edges, respectively (see \cite{SederbergZhengBakenovNasri2003} or \cite[Section 5.1]{Morgenstern2017} for details). They are used for the spline construction as follows: We denote the B-spline function $N_{1,3}^{\knots_u(z)}$ from \autoref{eq: cox de boor} above by $B_z^u$, applied with the knot vector $\knots_u(z)$, $n=1$ and $\degree=3$, and analogously by $B_z^v=N_{1,3}^{\knots_v(z)}$, the B-spline from \autoref{eq: cox de boor} applied with the knot vector $\knots_v(z)$.

\begin{figure}[ht]
 \centering
 \begin{tikzpicture}[scale=.8]
\node at (1.5,2.5) {\includegraphics[width=2.4cm]{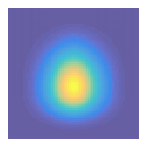}};
\draw[thick] (0,1) grid (4,5)
      (0,1.5)--(3,1.5) (1.5,1)--(1.5,2);
\draw[<->] (-.5,5.2) node [above] {$v$} |- (4.1,0) node [right] {$u$};
\fill[blue] (1.5,2) circle (.075);
\foreach \a/\b in {0/2,1/2,1.5/2,2/2,3/2,1.5/1,1.5/1.5,1.5/3,1.5/4}
\draw[thick,blue] (\a-.1,\b-.1)--++(.2,.2) (\a-.1,\b+.1)--++(.2,-.2);
\foreach \a/\b in {0/0,1/0,1.5/0,2/0,3/0,-.5/1,-.5/1.5,-.5/2,-.5/3,-.5/4}
\draw[thick,blue] (\a-.1,\b-.1)--++(.2,.2) (\a-.1,\b+.1)--++(.2,-.2);
\node[above] at (1.5,0) {\includegraphics[width=2.4cm]{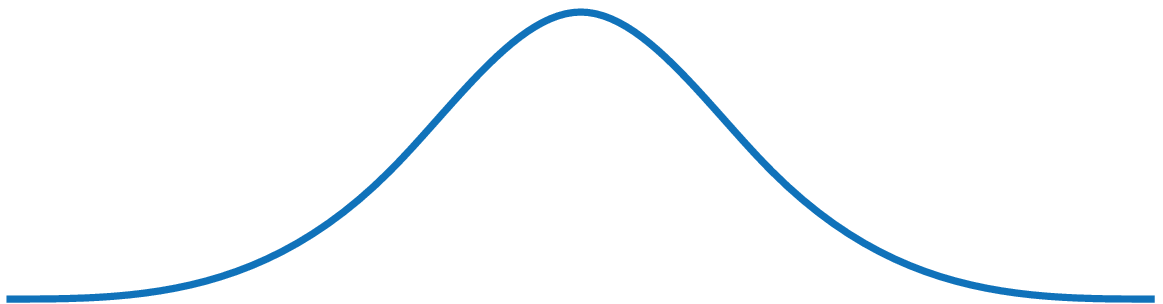}};
\node[left] at (-.5,2.5) {\includegraphics[height=2.4cm]{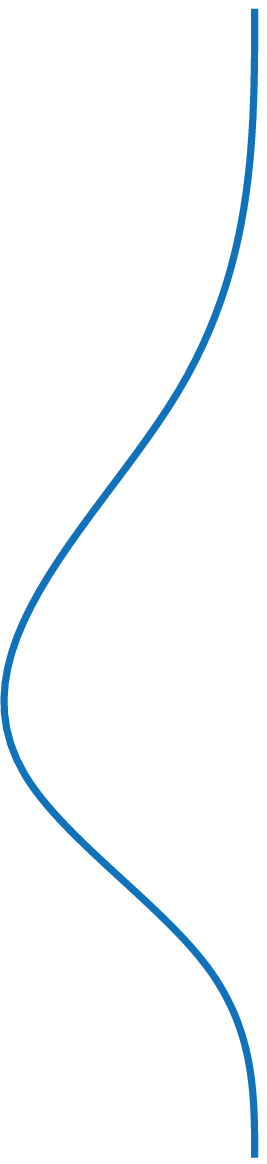}};
\end{tikzpicture}%
\caption{The T-spline function associated to the vertex is the product of univariate B-splines, constructed via distances to the nearest horizontal/vertical edges in a vertical/horizontal parameter direction.}
\label{fig: t-spline contruction}
\end{figure}
The spline function associated with the vertex $z$ is defined by $B_z(u,v)=\frac{w_z B_z^u(u)B_z^v(v)}{\sum_{r\in\nodes}w_r B_r^u(u)B_r^v(v)}$, also see \autoref{fig: t-spline contruction}.
\subsubsection{T-splines surfaces}
Each control point in the NURBS framework is associated with a pair of knot values coming from the two knot vectors. When a new knot is inserted, for example, for refinement purposes, unnecessary control points have to be added everywhere along a mesh line crossing the other direction to maintain the tensor product structure. Consequently, the control mesh is overrefined: this is a major drawback, particularly in the case of observations for which the risk of noise fitting is high. This may lead to oscillations in the fitted surface \cite{BraccoGiannelliGrossmannSestini2018} as the control points are determined by highly localized data points. In order to overcome this limitation, inserting a single control point without the need to add an entire row or column of control points should be made possible. The T-splines address this issue; they are described as point-based splines, i.e. a basis function of the spline space is defined for every vertex. Each of the basis functions comes from some tensor-product spline space. The position of a surface point on the T-splines surface is given by
\begin{equation}\label{eq: tspline surface}
 S(u,v) = \sum_{z\in\nodes}B_z(u,v)P_z,
\end{equation}
where $P_z$ and  $w_z$ represent the control point and weight associated with the vertex $z$. The weights are set to 1 in most applications with AST-splines, since the overlap of their supports is bounded a priori. Interested readers should refer to \cite{Casqueroetal2020} for an advanced generalization of AST-splines or \cite{Weietal2017} and \cite{Laietal2017} for some applications.
\subsubsection{Refinement strategy}\label{sec: mesh refinement}
We require that 
\begin{inlinelist}
\item the T-spline surface
\item the control points far away from the refinement location remain unchanged for the local refinement of the parametric mesh.  
 \end{inlinelist}%
 
\begin{figure}[ht]
 \centering
 \begin{tikzpicture}
 \pgfmathsetmacro{\c}{3.5}
\path (0,0) coordinate (1) ++(\c,0) coordinate (2) ++(\c,0) coordinate (3) (-.5*\c,-\c) coordinate (4) ++(\c,0) coordinate (5)
++(\c,0) coordinate (6) ++(\c,0) coordinate (7) (0,-2*\c) coordinate (8) ++(\c,0) coordinate (9) ++(\c,0) coordinate (10);
\foreach \a in {1,...,10}
\node[draw=gray, inner sep = 0pt,outer sep=3pt] (n\a) at (\a) {\includegraphics[width=3cm]{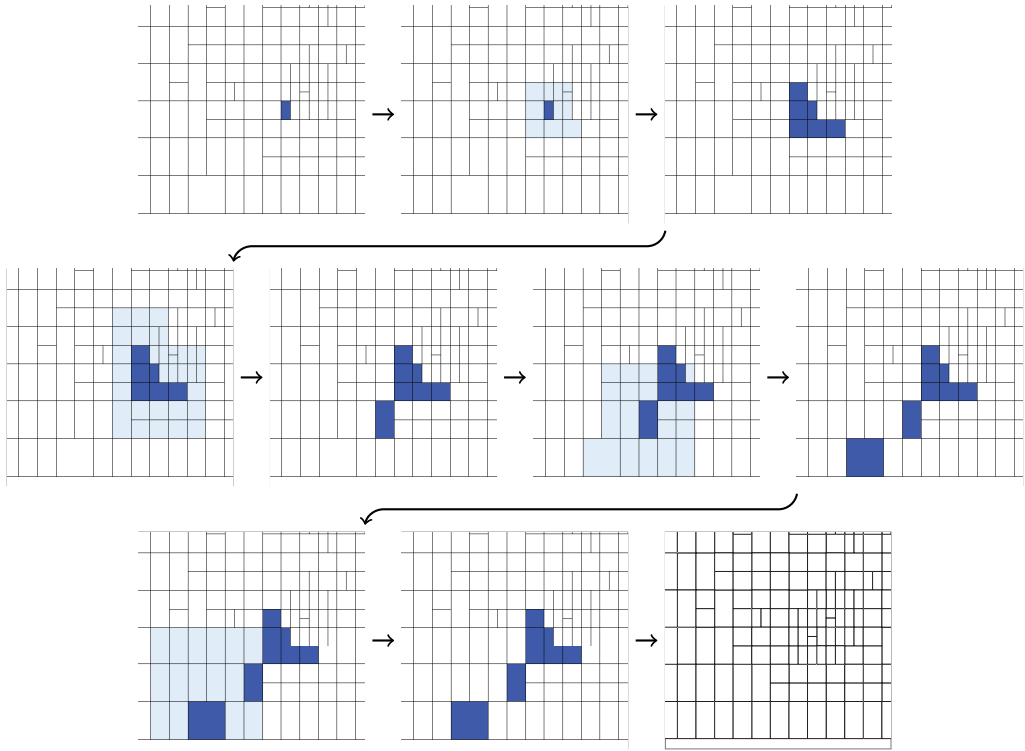}};
\foreach \a in {1,2,4,5,6,8,9}
{  \pgfmathsetmacro{\ap}{int(\a+1)}  \draw[thick,->] (n\a)  to (n\ap);  }
\foreach \a in {3,7}
\draw[thick,->] (n\a.south) ++(-1.5,0)  arc(-10:-90:.5*\c-1.5) -- ++(-1.5*\c,0) arc(90:170:.5*\c-1.5);
 \end{tikzpicture}
 \caption{Refinement scheme from \cite{MorgensternPeterseim2015}.
 Given a cell marked for refinement, the algorithm additionally marks all coarser cells in its neighborhood. The refinement routine is called recursively, checking their neighborhoods for even coarser cells, which are marked as well. The third iteration leaves the set of marked cells unchanged, and the algorithm ends with the refinement of all marked cells.
 The refinement direction for each cell depends on the refinement level.}
 \label{fig: 2D t-spline refinement}
\end{figure}

 We apply the refinement strategy from \cite{MorgensternPeterseim2015}, which is based on two key ideas. 
Firstly, cells with an even refinement level are refined via bisection with a vertical line, and cells with odd refinement level by a horizontal line. Thus, the refinement produces rectangles for odd refinement levels and squares for even refinement levels for the initial mesh consisting of squares. 
Secondly, the refinement of a cell requires that the cells in a certain neighbourhood be of at least the same refinement level (i.e. not coarser). The size of this neighbourhood depends on the polynomial degree, which is 3 in our setting. This second rule enforces a distance of $\approx p$ cells between T-junctions of different directions, and it is realized by a recursive refinement scheme that repeatedly checks the neighbourhood for coarser elements, see \autoref{fig: 2D t-spline refinement}. 
For details of the algorithm, we refer to \cite{MorgensternPeterseim2015}.

\subsection{Adaptive surface approximation}
We exploit the linear independence of the analysis-suitable T-splines to fit a surface to unevenly distributed points.
Our algorithm starts with an LS approximation of the data points, followed by an explicit local refinement strategy applied iteratively to cells in which the norm of the point-wise residuals exceeds a certain threshold.

We assume that we are given parametrized observations $\ell(u_j,v_j)\in\mathbb R^3$ for a set of parametric points $(u_j,v_j)\in D$, with $j=1,\dots,\nobs$, and a mesh with associated T-splines, as described in \autoref{sec: tsplines}.

\subsubsection{LS Principle}
The LS approach to surface approximation seeks control points $\{P_z\mid z\in\nodes\}$ such that the surface $S: D\to\mathbb R^3$ from \autoref{eq: tspline surface} minimizes
the sum of the squared euclidean distances between the observations $\ell(u_j,v_j)$ and their approximations $S(u_j,v_j)$, i.e.
\begin{equation}\label{eq: minimization LS}
 \min_{\{P_z\mid z\in\nodes\}\subset\mathbb R^3} \sum_{j=1}^\nobs \left\|S(u_j,v_j)-\ell(u_j,v_j)\right\|_2^2,
\end{equation}
with $\|\cdot\|_2$ as the euclidean distance. As the coordinates of any point on the T-spline surface are represented as the linear combination of the control points $P_z$, a perfect surface fitting would satisfy an overdetermined linear system $AP=L$, where $A$ is a $\nobs\times n$-matrix with entries $a_{ij}=B_j(u_i,v_i)$, $P\in\R^{n\times3}$ is a matrix containing the coordinates of the control points, and $L\in\R^{\nobs\times3}$ is the list of observations $\ell(u_i,v_i)$, with $\nobs$ being much larger than $n$. The LS solution is obtained by multiplying by $A\transp$ and solving the problem $\tilde AP=\tilde L$, where $\tilde A=A\transp A\in\R^{n\times n}$ is the so-called design matrix with entries 
$a_{jk} =\sum_{i=1}^\nobs B_j(u_i,v_i) B_k(u_i,v_i)$
and $\tilde L=A\transp L\in\R^{n\times3}$
is a matrix with lines $\ell_j=\sum_{i=1}^\nobs\ell(u_i,v_i)B_j(u_i,v_i)$. The LS solution $P$ obtained minimizes the residual $\left\|AP-L\right\|_2$.

Since the $i$-th line of $AP$ for any $i=1,\dots,\nobs$ reads 
$\sum_{j=1}^n B_j(u_i,v_i)\cdot P_j=S(u_i,v_i)$, this solves the minimization problem \eqref{eq: minimization LS}.

The refinement strategy is driven by the pointwise residual as an error indicator,
\begin{equation}\label{eq: error indicator}
e_i=\left\|S(u_i,v_i)-\ell(u_i,v_i)\right\|_2,\quad\text{for }i=1,\dots,\nobs.
\end{equation}
We note that a fairness function may be added to \eqref{eq: minimization LS} in order to obtain smooth surfaces; its use is extensively discussed in \cite{SkyttBarrowcloughDokken2015}. Data gaps and/or knot vector refinement for NURBS often lead to the occurrence of spline basis functions that have few or no observations in their support. Although the spline basis is linearly independent, the design matrix may be singular due to this challenging distribution of observations.  Ripples may appear in the surface approximation when the control points are determined by too few data points and/or too localized data points regarding the entire support. In order to overcome such situations, we make use of the Matlab function \texttt{mldivide}\footnote{see also \url{https://mathworks.com/help/matlab/ref/mldivide.html}}, commonly known as Matlab's backslash operator, which automatically chooses an appropriate method to solve the problem using a QR, Cholesky, triangular, Hessenberg or LU solver, to cite but a few. Unfortunately, challenging cases still arise and are accentuated as the number of iterations increases, leading to numerical inaccuracies or slow convergence. Applying LS approximation in the first step for a high quality of reconstruction and switching to a local explicit method was shown to be a good alternative by producing well-behaved surfaces in the case of outliers, gappy, scattered and noisy data, see, for example, \cite{BraccoGiannelliGrossmannSestini2018}. As these are the main characteristics of low-cost laser scanners for land remote sensing applications, this approach combines 
  the power and flexibility of T-splines, which are widely used in the CAD context, utilizing 
  a local explicit method for data approximation with challenging density, as arises within a GIS context. This procedure 
  speeds up the fitting without greatly affecting the quality of the approximation. 
In the next section, we will describe the MBA principles in more detail.
\subsubsection{MBA principle}
The multilevel B-spline approximation was introduced in \cite{LeeWolbergShin1997} for data interpolation and extended by \cite{ZhangTangLi1998} for approximating scattered data with HB-splines. 
By construction, the MBA technique is particularly appropriate when data gaps occur, as shown extensively in \cite{SkyttBarrowcloughDokken2015}.

\begin{figure}[!ht]
\centering
\includegraphics[height=.7\textheight]{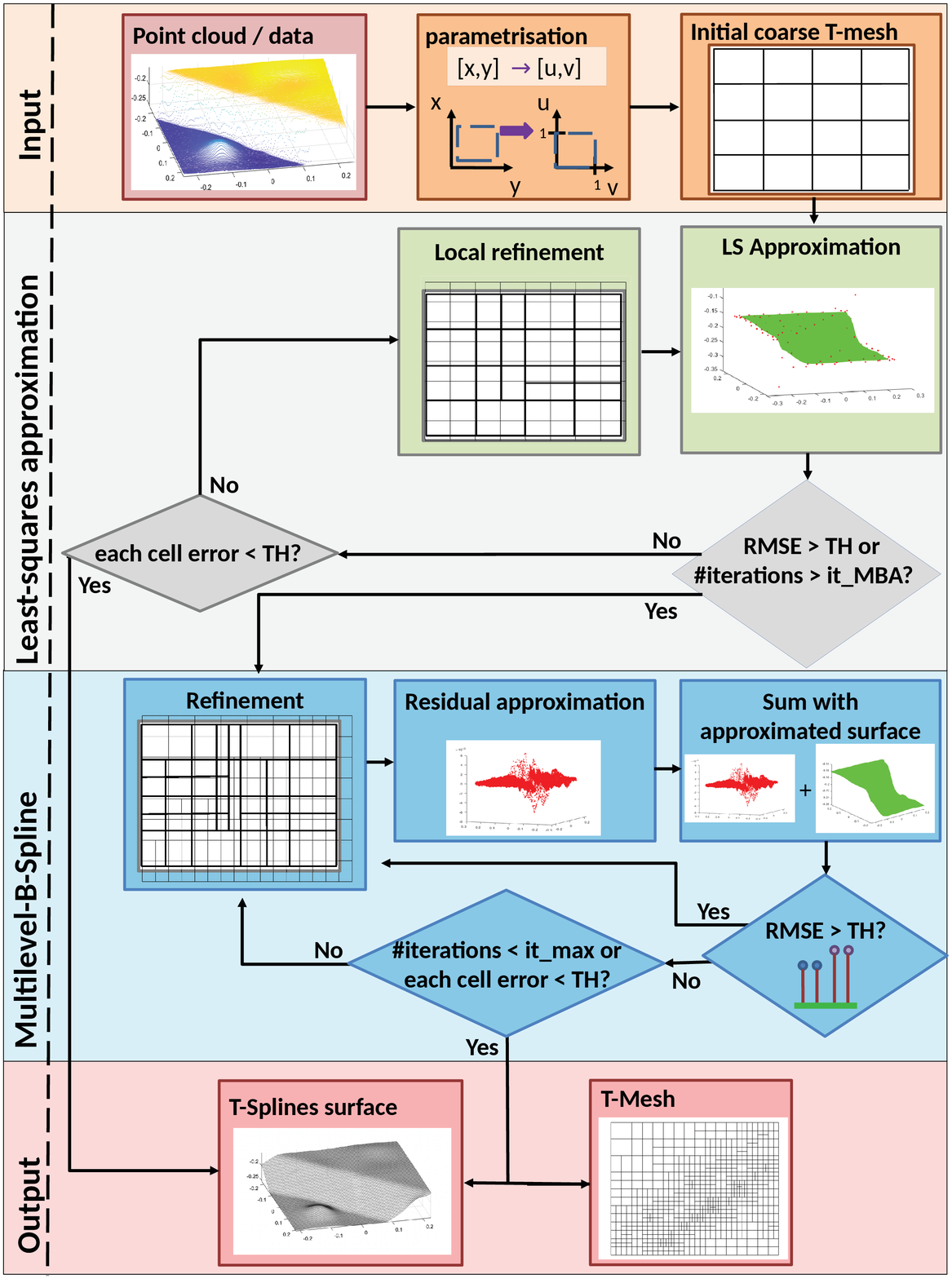}
 \caption{Combination of the MTA and LS-T for adaptive surface approximation with T-splines.}
 \label{fig: flow chart}
\end{figure}

Here, we apply it with T-splines. The MTA can be summarized as follows (see also \autoref{fig: flow chart}): 
we assume that all observations are of the form $\ell(u,v)=\bigl(u,v,z_\ell(u,v)\bigr)$, and that we have a preceding surface approximation $S:D\to\R^3$ via the MTA or LS method, which is also of the form $S(u,v)=\bigl(u,v,z_S(u,v)\bigr)$ (see \autoref{fig: flow chart}, top row). 
A mesh refinement is performed, driven by the error indicators $e_i$ from the \autoref{eq: error indicator} above, yielding a new set of spline functions $\bigl\{B_i\mid i\in\{1,\dots,N\}\bigr\}$.
Instead of computing the control points of the new surface approximation using LS (\autoref{fig: flow chart}, second part ``Least-squares approximation''), they are calculated individually using the residuals of the preceding approximation (\autoref{fig: flow chart}, third part ``Multilevel B-spline''). We call $z_c= z_S(u_c,v_c)-z_\ell(u_c,v_c)$ the residuals at observation points $(u_c,v_c)$, $c=1,\dots,\nobs$. 
The coefficients $q_i$ of the approximated residual surface, involving a given threshold $\thr$ for the approximation, are determined by 
\begin{align}\label{eq: MBA coefficients phi}
 \phi_{i,c} &= \frac{B_i(u_c,v_c)z_c}{\sum_{j=1}^NB_j(u_c,v_c)^2}
 \quad\text{for all $i=1,\dots,N$ and $c=1,\dots,\nobs$},\\[1ex]
 \label{eq: MBA coefficients q}
 q_i&=\begin{cases}
0       & \text{if }\left|S(u_c,v_c)-z_c\right|<\thr\\&\text{for all }(u_c,v_c)\in\supp B_i,\\
\displaystyle\frac{\sum_{c=1}^\nobs B_i(u_c,v_c)^2\phi_{i,c}}{\sum_{c=1}^\nobs B_i(u_c,v_c)^2} & \text{otherwise}
      \end{cases}
 \quad\text{for all $i=1,\dots,N$.}
\end{align}
Both \autoref{eq: MBA coefficients phi} and \autoref{eq: MBA coefficients q} are themselves results of localized LS minimizations:  
\autoref{eq: MBA coefficients phi} minimizes $\sum_{i=1}^N\phi_{i,c}^2$ provided that $\sum_{i=1}^N\phi_{i,c}B_i(u_c,v_c)=z_c$, which can be seen as a localized Ridge regression%
, and
\autoref{eq: MBA coefficients q} minimizes $\sum_{c=1}^\nobs\bigl(q_iB_i(u_c,v_c)-\phi_{i,c}B_i(u_c,v_c)\bigr)$.
In practise, the sums are only computed for those indices $(i,c)$ where $(u_c,v_c)$ is in the support of $B_i$, since all other contributions are zero.
The new T-spline surface $S_1:D\to\R^3$ is computed as the sum $S_1(u,v)=\bigl(u,v,\,z_S(u,v) {+} R(u,v)\,\bigr)$ of the initial surface $S$ and the surface obtained from the residual approximation $R:D\to\R$, $R(u,v)=\sum_{i=1}^nB_i(u,v)\cdot q_i$. The residual approximation is added to the surface approximation incrementally at each step, until the maximum number of iterations is reached (\autoref{fig: flow chart}, last part ``Output'').

Note that the first case in \autoref{eq: MBA coefficients q} also applies if the support of $B_i$ does not contain any data point $(u_c,v_c)$.
Hence, fine-scale coefficients are, by construction, set to zero where data points are missing or have already been met sufficiently well by a coarse-scale approximation, which also follows the spirit of sparse sensing, where superfluous coefficients are set to zero by minimizing the L2 norm (Ridge regression) or L1 norm (LASSO Algorithm) of the coefficient vector. 
The approximated residual surface is, thus, sparse and the computation time decreases significantly compared to the LS strategy.
The adaptive choice of fine-scale basis functions is combined with the adaptive mesh refinement explained in \autoref{sec: mesh refinement} to reduce the generation of data structures for which the coefficients are set to zero in the approximation step.

\subsubsection{Iterative surface approximation using LS or MBA}
The starting point of the iterative surface approximation is a tensor product B-spline surface defined over a coarse mesh. An LS adjustment is performed iteratively in the element regions, where one observation is associated with an error term higher than a given threshold $\thr$. Intuitively, few cells are refined with a large threshold. On the contrary, all cells will be halved for a very low $\thr$, which results in a global refinement, similar to that for NURBS with a dyadic mesh structure. 

The procedure of scattered data approximation is summarized in \autoref{alg: fitting}.
The choice of $\thr$ is left to the user's convenience, depending on the level of accuracy needed, the computation time that can be supported or the need to avoid noise fitting. 
$\thr$ can be fixed in order not to exceed the variance of the observations in the $z$-direction. It should be chosen carefully with the aim of avoiding noise fitting.
Here, the approximation will cease to improve when increasing the level of refinement, which indicates an optimal choice of $\thr$. Two loops are visible in \autoref{alg: fitting}, both containing a mesh refinement. Hence throughout this paper, refinement is always implied when we mention an iteration of the algorithm.
Since we do not use any iterative schemes in sub-steps of the algorithm, we refer to  the refinement degree increases with each iteration.

We also mention an important point: we cannot and will never expect the MBA to perform as well as an LS method. However, the MTA avoids the drawback of a global LS adjustment as the iteration level increases. Additionally, the approximated surfaces are well-behaved, which is particularly interesting for land remote sensing applications. None of the surface approximation techniques will reproduce the surface exactly as the truth does not exist for real applications. 

\begin{algorithm}
 \begin{algorithmic}
  \While{ there exist points $i$ with $e_i>\thr$ \& max.\ number of iteration not reached}
  \State{refine the surface}
  \State{compute the approximate surface (LS or MBA)}
  \State{compute the quality parameter (residual, $\rmse$, points outside tolerance, error term)}
  \EndWhile
 \end{algorithmic}
\caption{Iterative algorithm for T-splines surface approximation with LS or MBA.}
\label{alg: fitting}
\end{algorithm}

\section{Simulations}\label{sec: simulations}
In this section, we validate the use of a combination of LS and MTA for different simulated point clouds. Accordingly, we generate noisy and scattered surfaces with and without sharp edges, outliers and data gaps. We further define the criteria by which the fitted surface will be judged.
\subsection{Simulated point clouds}\label{sec: point clouds}
\subsubsection{Reference surfaces}
The reference T-spline surfaces correspond to smooth geometries given by 
\begin{align}
\Ssmooth(u,v) &= s_1^{(9)}(u,v) + s_2(u,v) + s_3(u,v),\\[1ex]
\text{and}\quad
\Ssharp(u,v) &= s_1^{(30)}(u,v) + s_2(u,v) + s_3(u,v),\\[1ex]
\text{with}\quad
s_1^{(k)}(u,v) &= \tfrac16\bigl(\tanh(k(v-u)) + 1\bigr), \\
s_2(u,v) &= 0.1 \exp\Bigl(-30\bigl((u-0.415)^2+(v+0.415)^2\bigr)\Bigr), \\
\text{and}\quad
s_3(u,v) &= 
\begin{minipage}[t]{.8\textwidth}
\scriptsize $ - 0.03\exp\Bigl(-20\bigl((u+0.5)^2+(v-0.5)^2\bigr)\Bigr)  + 0.03\exp\Bigl(-10\bigl((u+0.6)^2+(v-0.6)^2\bigr)\Bigr) \\
    - 0.03\exp\Bigl(-10\bigl((u+0.4)^2+(v-0.6)^2\bigr)\Bigr) 
    + 0.02\exp\Bigl(-10\bigl((u+0.6)^2+(v-0.4)^2\bigr)\Bigr) \\
    + 0.01\exp\Bigl(-10\bigl((u+0.7)^2+(v-0.3)^2\bigr)\Bigr) 
    + 0.02\exp\Bigl(-10\bigl((u+0.1)^2+(v-0.7)^2\bigr)\Bigr) \\
    - 0.01\exp\Bigl(-20\bigl((u+0.6)^2+ v    ^2\bigr)\Bigr)$
\end{minipage}
\end{align}


\begin{figure}[ht]
\centering
\includegraphics[height=.21\textheight]{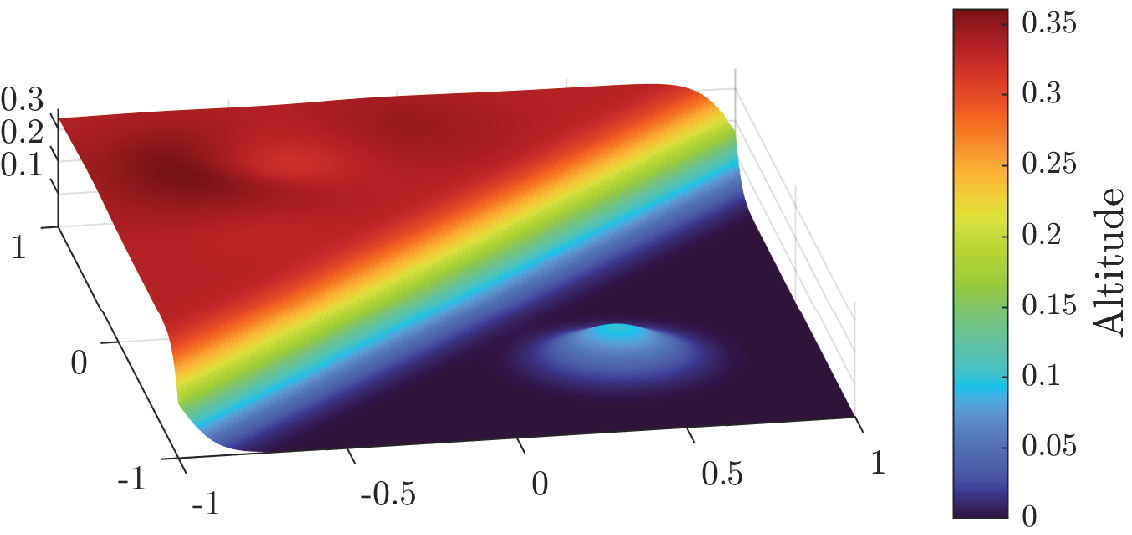}~%
\includegraphics[height=.21\textheight]{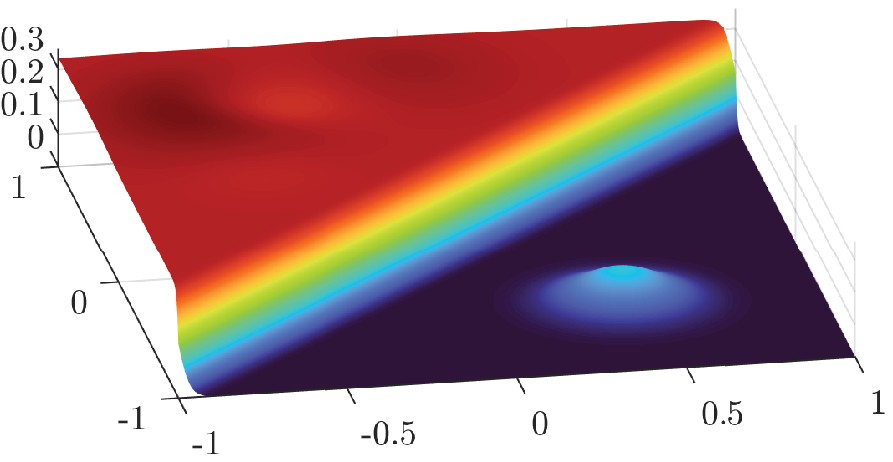}%
\caption{Surface $\Ssmooth$ and $\Ssharp$ serving as references.}
\label{fig: reference surface}
\end{figure}

The surface consists of three components: 
$s_1^{(k)}(u,v)$ corresponds to a dam, which can be further modified by varying the coefficient $k$ of the $\tanh$ argument to challenge the refinement algorithm by, for example, a steeper descent.
$s_2(u,v)$ corresponds to a Gaussian bell, i.e. a hill in real life.
The smooth ripples $s_3(u,v)$ are added to avoid a simple plane in the upper part of the surface and were favourable for the 3D printing of this surface, as presented in \autoref{sec: real data}.
The range of values for $u$ and $v$ is between $-1$ and $1$ for the generation of the surface. We sample both surfaces with $\nobs=200\times200=\num{40000}$ data points, in order to compare the convergence behavior of the proposed approximation schemes. 

\paragraph{Data gaps and outliers}
In order to investigate the extent to which data gaps affect the approximation, we selected a patch on $\Ssmooth$ given by $[-\frac14,0]^2$ for which the corresponding observations were set to 0. This leads to a small square subdomain on the dam where data is missing. The surface will be called $\Sgaps$. We will further introduce 5\%  of outliers in the $z$-component added randomly and call the surface obtained $\Soutliers$. The outliers follow a student distribution and can reach up to 10 times the maximum $z$-component of the original point cloud.
\subsubsection{Surface approximation}
All four surfaces generated are noised, obtaining noised surfaces $\tSsmooth$,$\tSsharp$,$\tSgaps$,$\tSoutliers$. In the $z$-component, we add a Gaussian noise of a standard deviation of \num{0.003} generated with the Matlab function \texttt{randn}. The same procedure is applied to the $x$- and $y$-components with a Gaussian noise of a standard deviation of 0.001, i.e. we assume that the $z$-component in a real case would be noisier than the horizontal components.
Consequently, the observations generated are not exactly gridded in the $(x,y)$-plane and mimic real observations. We placed ourselves in the framework of Monte Carlo (MC) simulations and generated 100 different noise vectors. The surfaces are fitted using the LS-T and MTA methods proposed in \autoref{sec: maths}, which we summarize as follows:
\begin{itemize}
\item The LS-T only for all iterations
\item The LS-T for the first three iterations followed by the MTA
\end{itemize}
We chose a threshold $\thr=0.01$ to avoid overfitting. This threshold corresponds approximately to three times the standard deviation of the noise in the $z$-component. Cells are refined if they contain more than one point outside the tolerance. The maximum number of iterations is fixed to eight for LS-T and ten for MTA for computational reasons. Increasing the number of iterations does not change the conclusions of our work, but is extremely computationally demanding for the LS-T strategy using Matlab, as shown in \autoref{sec: results}. An implementation in C++ of the algorithms is to be performed in the near future.
The mean and standard deviation of the performance indicators defined in the next section are computed for each MC simulation and averaged over the 100 runs to judge the goodness of fit of the surface fitting.

\subsection{Performance indicator}\label{sec: performance indicator}
Judging the goodness of fit of a surface approximation is far from easy; it should contain both heuristic and statistic components. As a rough idea, we could require that the surface should be smooth without ``ripples'' with a maximum error term that does not exceed twice the threshold in the absence of outliers.  
In this contribution, we follow \cite{Wang2009} and define an optimal surface approximation as being 
\begin{inlinelist}
\item accurate such that the geometric error between approximated surface and point cloud is small
\item continuous even in case of sharp edges or data gaps
\item fair, i.e. without fluctuations coming from noise overfitting
\item concise with the right number of parameters
\end{inlinelist}.
Following \cite{SkyttHarphamDokkenDahl2017} or \cite{BraccoGiannelliGrossmannSestini2018}, we define following performance indicators: 
\begin{itemize}
\item The root mean square error between the noised surface and the approximation, defined as 
$\rmsenoisek=\frac1{\sqrt\nobs}\bigl\|\hat S_k-\tilde S\bigr\|_2$, where $\hat S_k$ is the approximate surface obtained after the $k$-th iteration,
\item the root mean square error referring to the exact surface without noise,\\ $\rmsemathk = \frac1{\sqrt\nobs}\bigl\|\hat S_k- S\bigr\|_2$,
\item the maximum error defined as $\maxerr[k] = \max\bigl\|\hat S_k-\tilde S\bigr\|_2$
\item the number of points outside tolerance $n_{\text{out},k}$,
\item the number of control points $\ncpk$ estimated for a given iteration $k$ of the refinement,
\item the computation time $\ct_k$. We used an 64-bit operating system with \SI{8}{\giga\byte} RAM and an Intel(R) Core(TM) i5-63000U CPU @ \SI{2.4}{\GHz}, \SI{5.8}{\GHz} 
\end{itemize}
These criteria were specifically chosen for comparing NURBS, MTA and LS-T approximation fairly.

\subsection{Results}\label{sec: results}
The results are presented with the aim of being both visual and statistical. Consequently, we present 
\begin{inlinelist}
 \item the T-meshes and fitted surfaces, highlighting heuristically the difficulties raised by the surface approximation 
 \item the means of the performance indicators obtained after 100 MC simulations computed for all four generated surfaces under consideration. 
\end{inlinelist}
The results are presented both in graphical and tabular forms. The variances of the performance indicators were below \num{1e-5} for the $\rmse$ and \num{1e-4} for $\maxerr$; they are not specified  so as not to overload this article. 

\subsubsection{Smooth surface}
\begin{figure}[ht]
 \centering 
 \includegraphics[width=.3\textwidth]{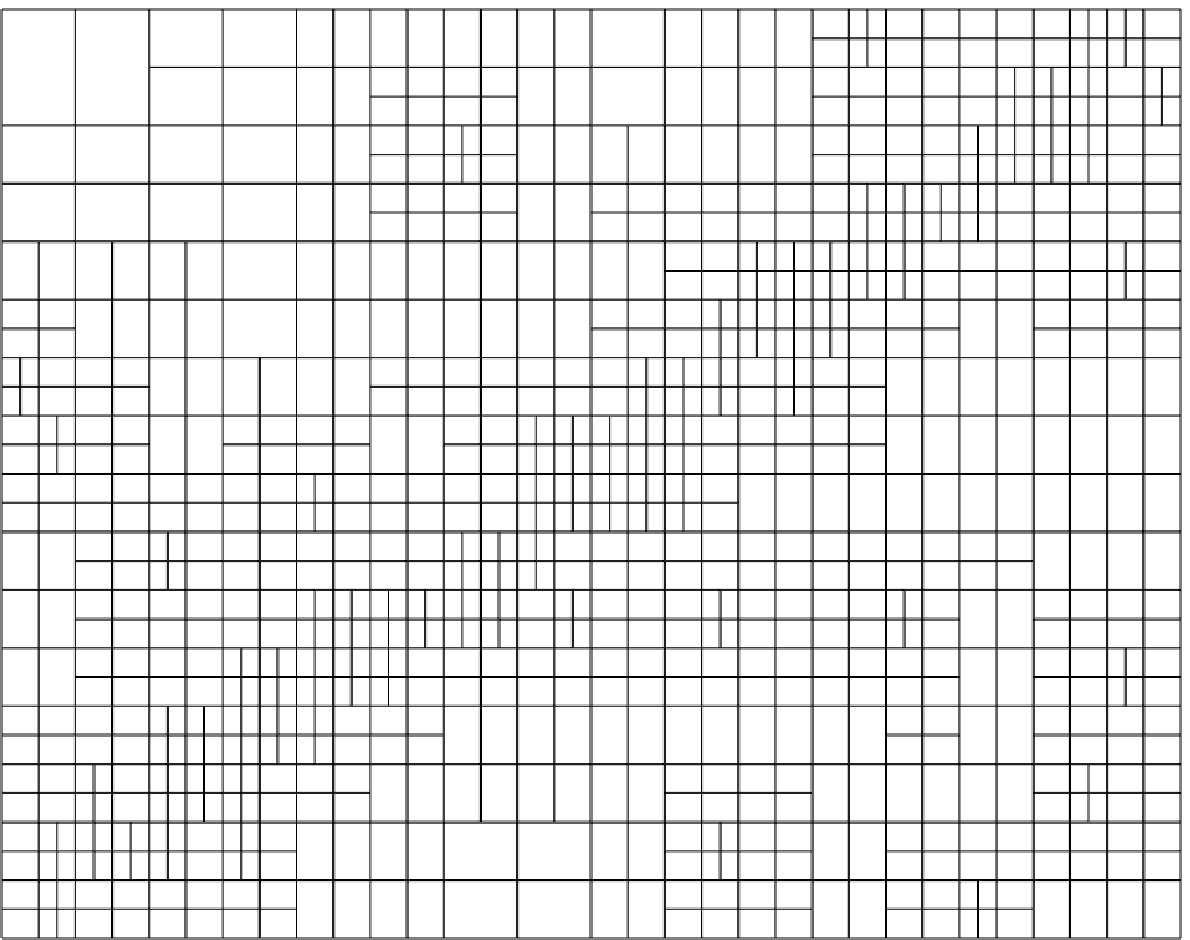}%
 \hspace{.049\textwidth}%
 \includegraphics[width=.3\textwidth]{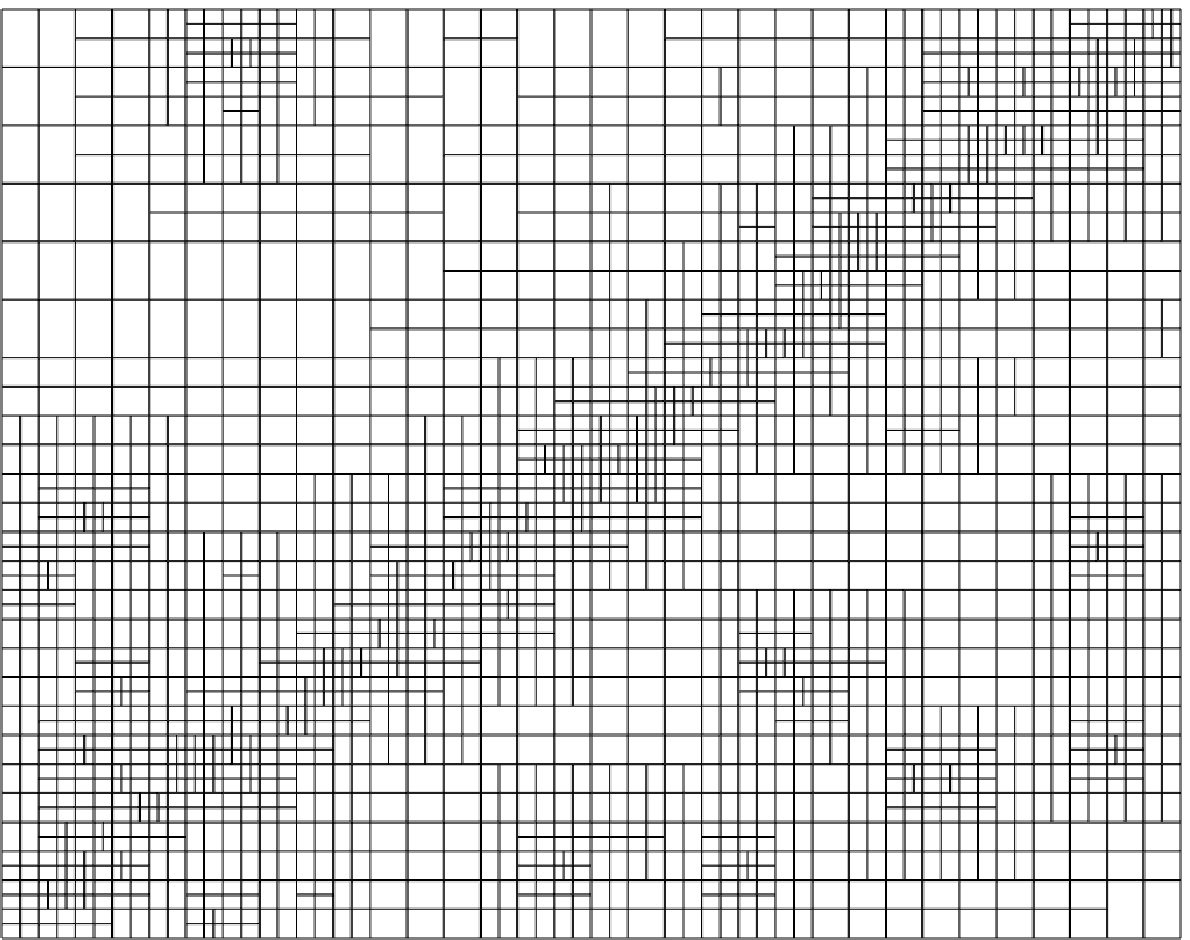}%
 \hspace{.049\textwidth}%
 \includegraphics[width=.3\textwidth]{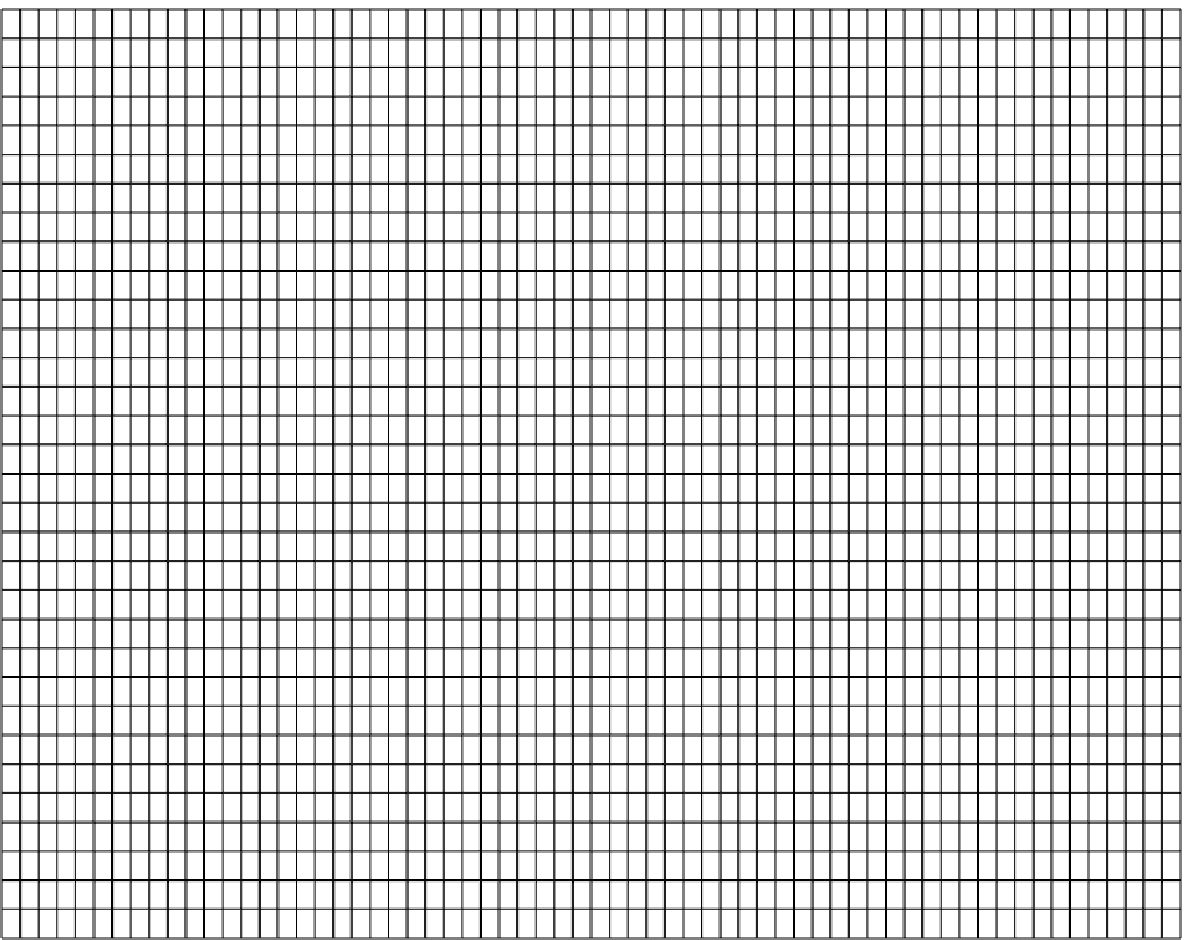}
 \caption{Approximations of $\Ssmooth$, T-meshes. Left: LS-T, middle: MTA and right: NURBS approaches.}
 \label{fig: smooth surface}
\end{figure}

We draw the following comments from \autoref{tab: approximation results} and \autoref{fig: smooth surface}, which presents the results of the performance indicators and the T-mesh for the LS-T, MTA and NURBS approaches, respectively:
\paragraph{T-meshes and fitted surface}
The number of maximum iterations for the MBA and LS/NURBS strategies were fixed to ten and eight, respectively, with the aim of performing a fair comparison between the two methods. 
The performance indicators (except for the number of control points) are similar for both methods, as \autoref{tab: approximation results} highlights.
Correspondingly, the number of control points estimated using the MTA is slightly higher than for the LS-T method (\autoref{tab: approximation results}). 
Fortunately, the computation time after ten iterations is similar and even smaller for the MTA regarding the LS-T strategy with eight iterations: 
this is easily explained as no matrix inversion has to be performed with the MBA approach.

The consequence of the increase of control points is an MTA-T-mesh with a finer structure compared to the LS-T-mesh for a similar $\rmse$. 
This is clearly visible in \autoref{fig: smooth surface} (left and middle): more lines are introduced, particularly along the dam, where a higher degree of refinement is needed to fit the surface properly. 
For the sake of completeness and fair comparison with usual T-splines fitting methods, we computed the $\rmse$ for the LS-T strategy with the same number of control points as for the MTA, see \autoref{tab: approximation results}. Besides the fact that the computation time strongly increases, the performance indicator stagnated, following the results of \cite{SkyttBarrowcloughDokken2015}. 
We performed the T-splines surface fitting with a global threshold of \num{0.01}. Decreasing the $\thr$ would have led to a higher level of refinement along the dam as no refinement is performed in the domain where the error term had already reached the threshold. 

We point out that the simulated Gaussian bell is not high enough to be seen in the refined mesh due to the added noise; a decomposition of the surface into patches to adapt the threshold locally, as proposed in \cite{FengTaguchi2017}, could lead potentially to an overfitting. This is avoided here by a global threshold and a non-over-refined mesh.  

\begin{figure}[ht]
\centering
\includegraphics[width=.45\textwidth]{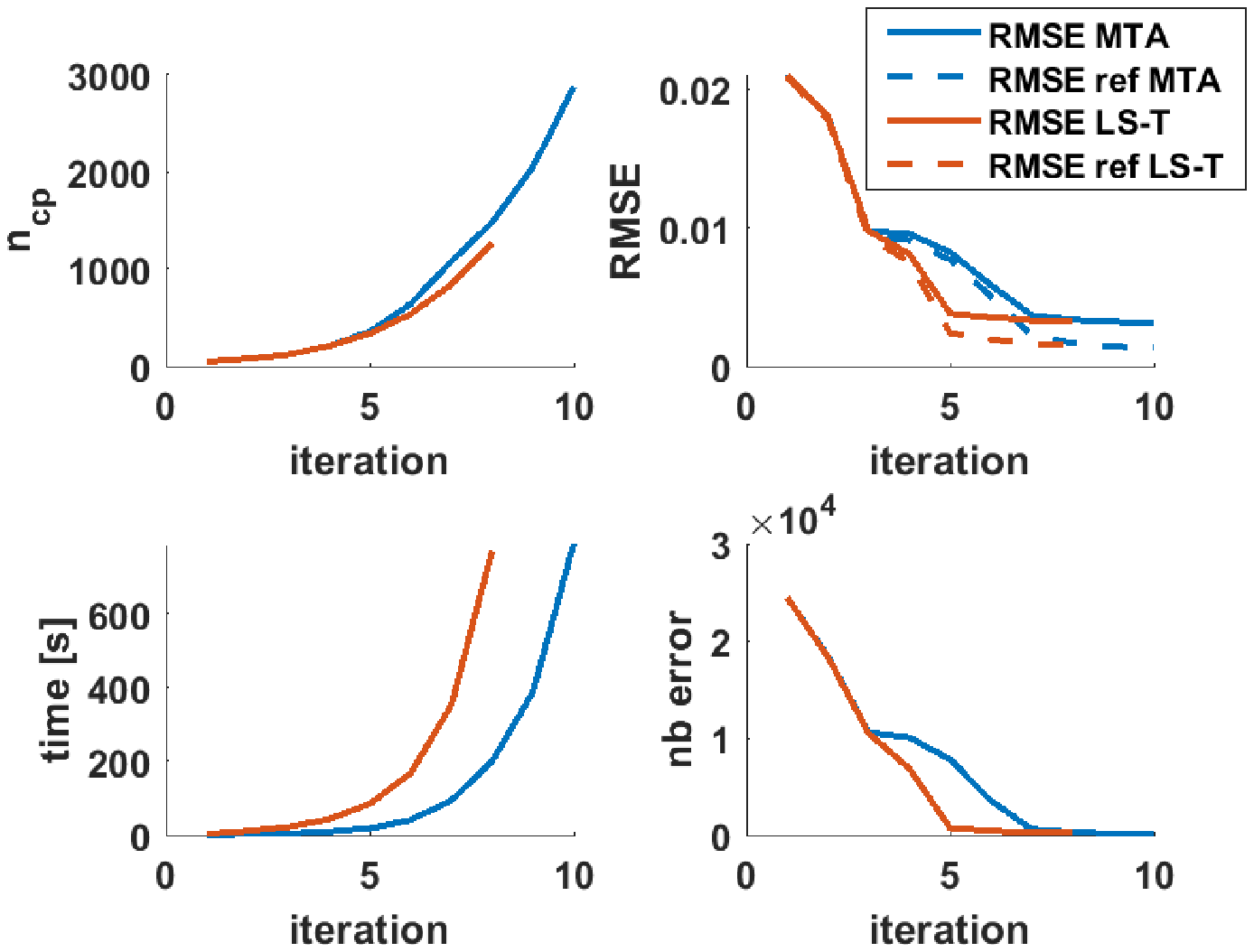}
\includegraphics[width=.45\textwidth]{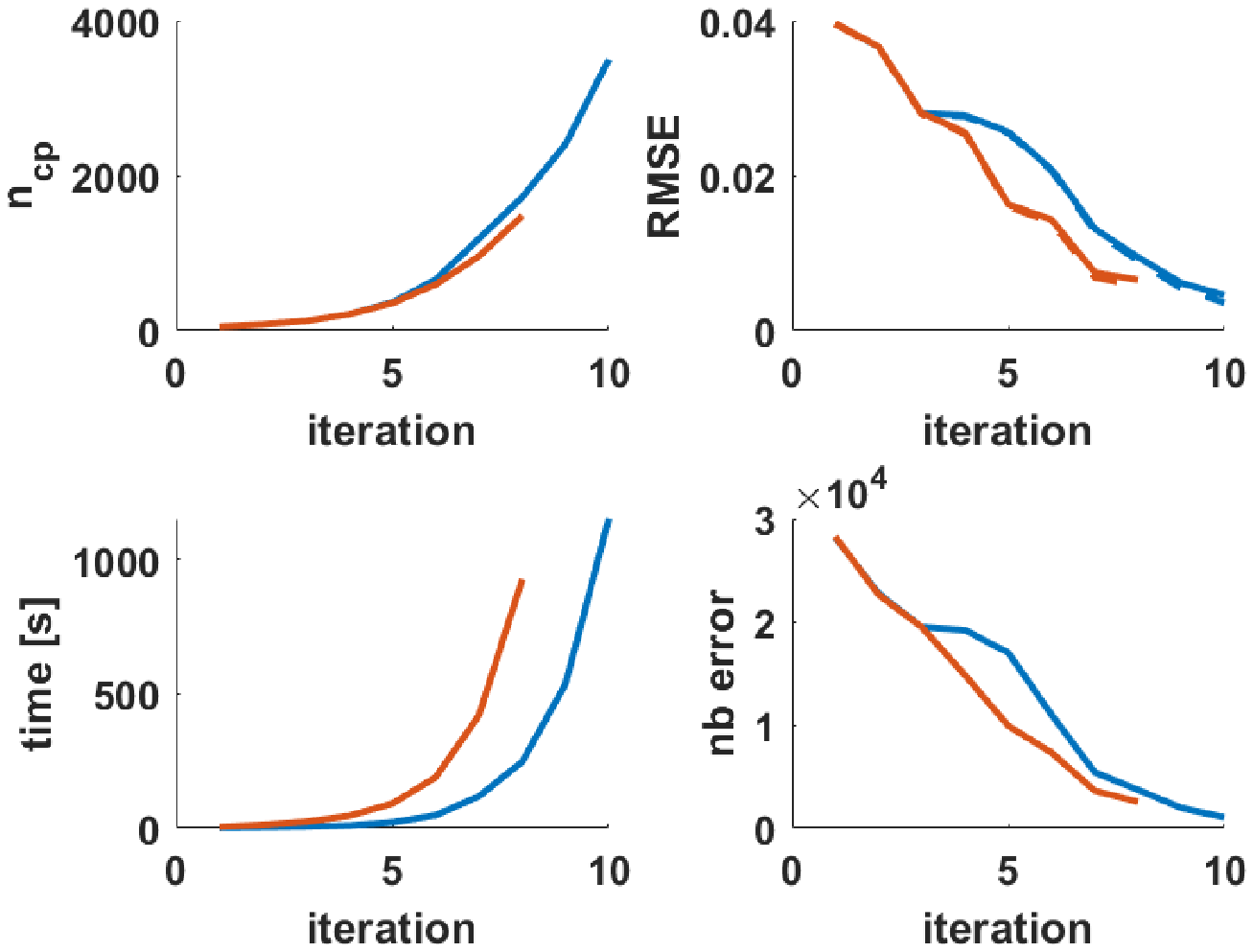}\\
\includegraphics[width=.45\textwidth]{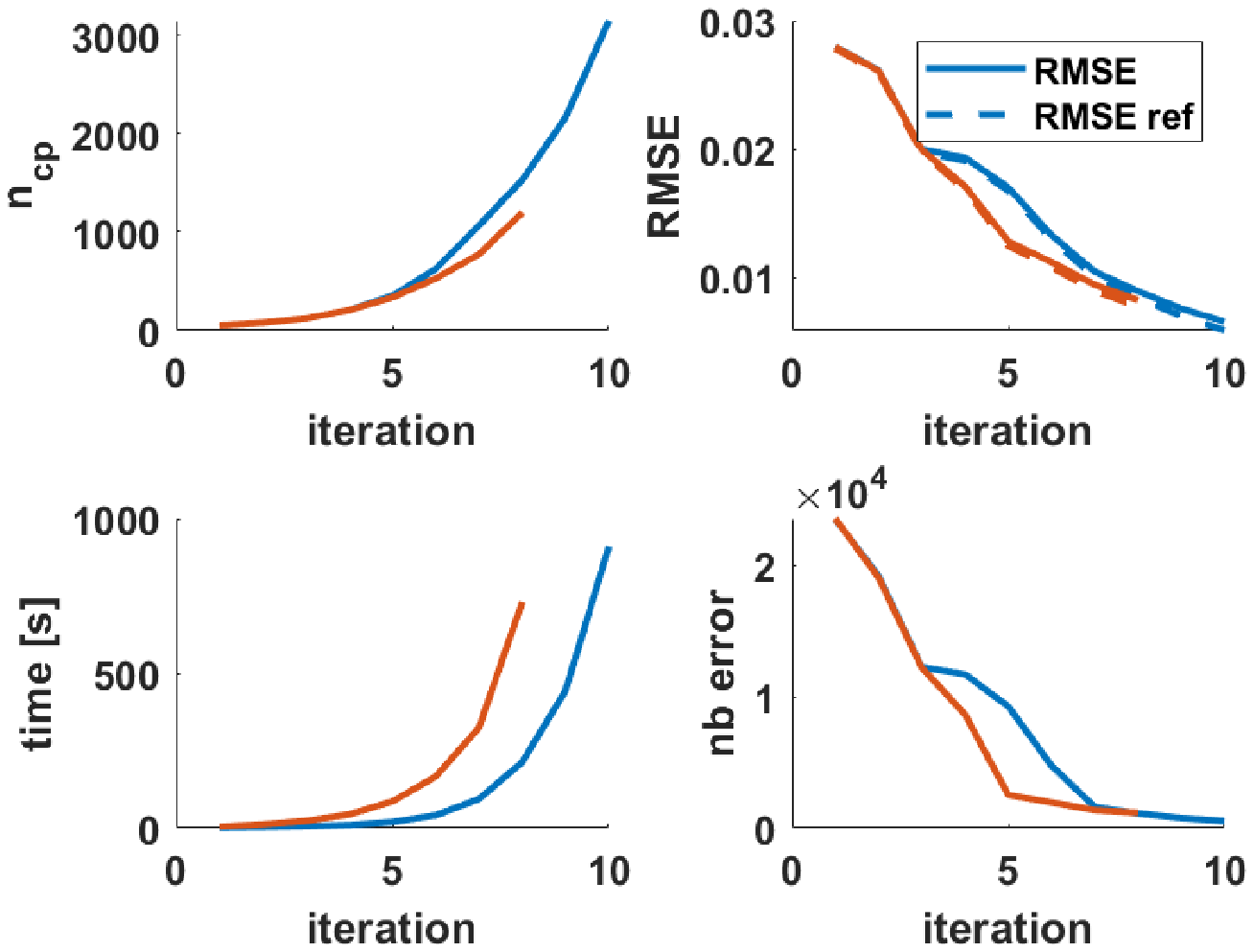}
\includegraphics[width=.45\textwidth]{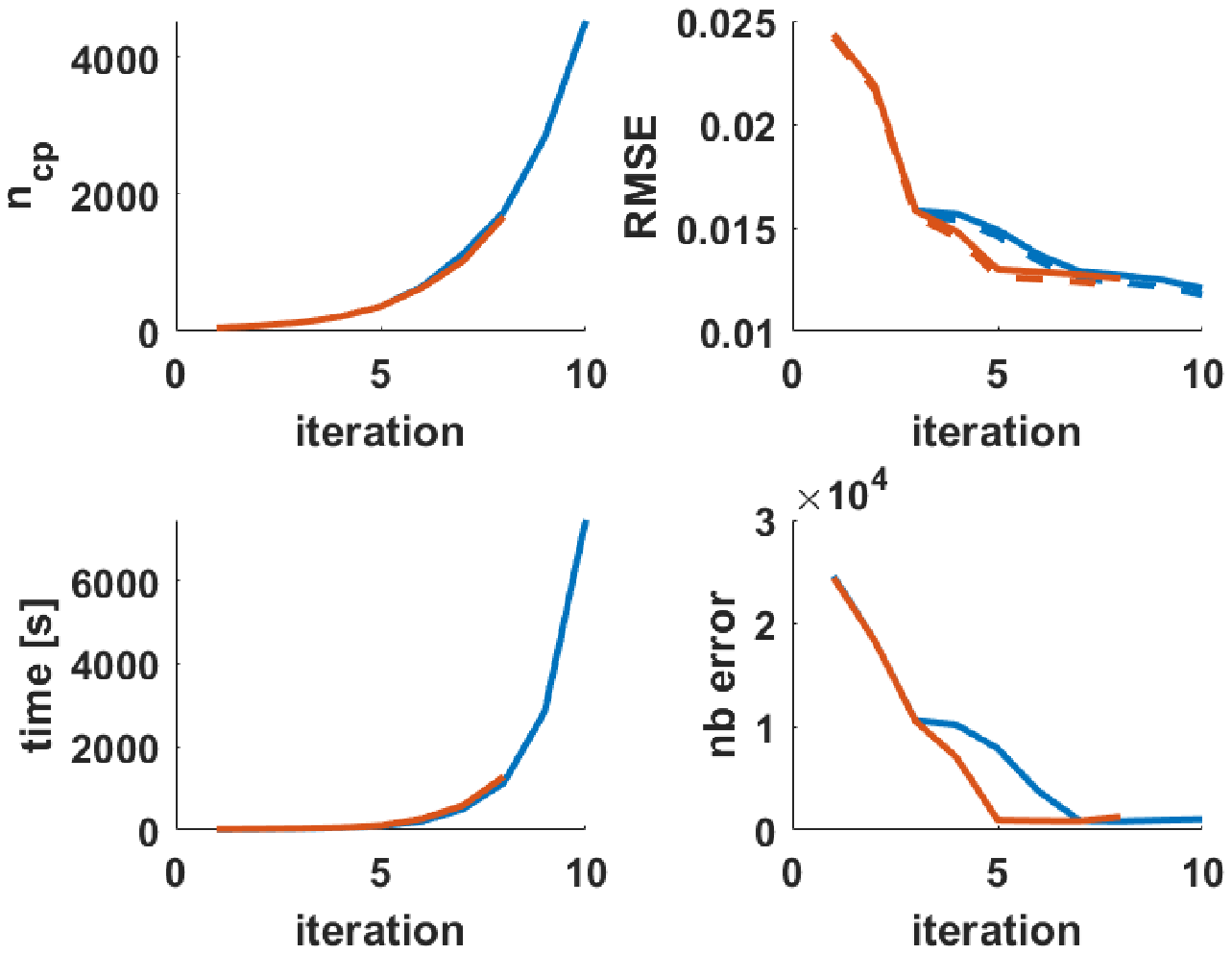}
\caption{Main performance indicator specified in \autoref{tab: approximation results} versus the iteration number. The blue curve corresponds to the MTA and the red to the LS-T. Top left: $\Ssmooth$, top right: $\Ssharp$, bottom left: $\Sgaps$ and bottom right: $\Soutliers$}
\label{fig: performance indicator}
\end{figure}

\paragraph{Performance indicator}
We mentioned above that for a given iteration step, $\ncpk$ is higher for MTA than for LS-T, which is clearly visible in the T-meshes. Correspondingly, the $\rmsemathk$  or $\rmsenoisek$ noise decreases with the number of iterations more rapidly for LS, as they follow $\ncpk$ versus the iteration step, see \autoref{tab: approximation results} and \autoref{fig: performance indicator}. 

We mention that both $\rmse$s decrease similarly and saturate to a common value close to the standard deviation of the simulated noise. 
Combined with a small number of points outside the tolerance regarding the total number of surface points (128 and 243 for the MTA and LS-T, respectively), these results highlight the goodness of the surface fitting, which is independent of the approximation strategy for a smooth surface. 

The $\rmse$ still decreases from a given number of iterations (\num 5 for LS-T and \num 7 for MTA, as shown in \autoref{fig: performance indicator}), but not as quickly as in the initial steps. 
This has to be linked to the exponential increase of the computation time; the latter should be weighed with the decrease of the $\rmse$. For the same number of control points to estimate, LS-T needs a much higher computation time than MTA.
Information criteria could be an answer to the challenge of determining the optimal number of iterations and the threshold \cite{KermarrecAlkhatib2020}. 
This remains the topic of a later contribution.
\paragraph{Comparison with NURBS}
In \autoref{sec: maths}, we highlighted that the NURBS is a particular case of the T-splines. The NURBS mesh is drawn in \autoref{fig: smooth surface} (right). It highlights the drawback of the global approximation, i.e. all cells are halved at each iteration resulting in a fine mesh as no local refinement is performed. 
Ripples and oscillations are more likely to arise in the NURBS surface due to overfitting. 
This is an unwanted effect, which can be avoided with a local refinement strategy. In addition, the computation time is higher, which is due to the number of control points which need to be estimated by LS. The preconditioned conjugate gradient method used by the Matlab function mdivide to speed up the computation reaches its limit here. 
The maximum error is twice higher, which could be due to the noise being fitted in the flat domain.

We mention that the surface could be smoothed with the help of a fairness functional \cite{Wang2009,ZhengWangSeah2005}.
The use of such additional functions necessitates the tuning or estimation of the fairness parameters, which increases the computation time significantly. 
Additionally, smoothing a surface in a real case may be unwanted. 
We do not follow this direction in this contribution. However, we are aware that this is a powerful option to obtain smooth fitted surfaces from noisy point clouds.
\subsubsection{Sharp surface}
We do not perform further investigations with NURBS in the following sections so as not to overload this contribution.
This strategy was shown to be suboptimal compared to local refinement strategies with T-splines in terms of the performance indicator and overfitted for a smooth and easy geometry. 
We will focus in this section on the comparison between LS-T and MTA for more challenging surfaces with sharp geometries and data gaps.
\begin{figure}[ht]
\centering 
\includegraphics[width=.3\textwidth]{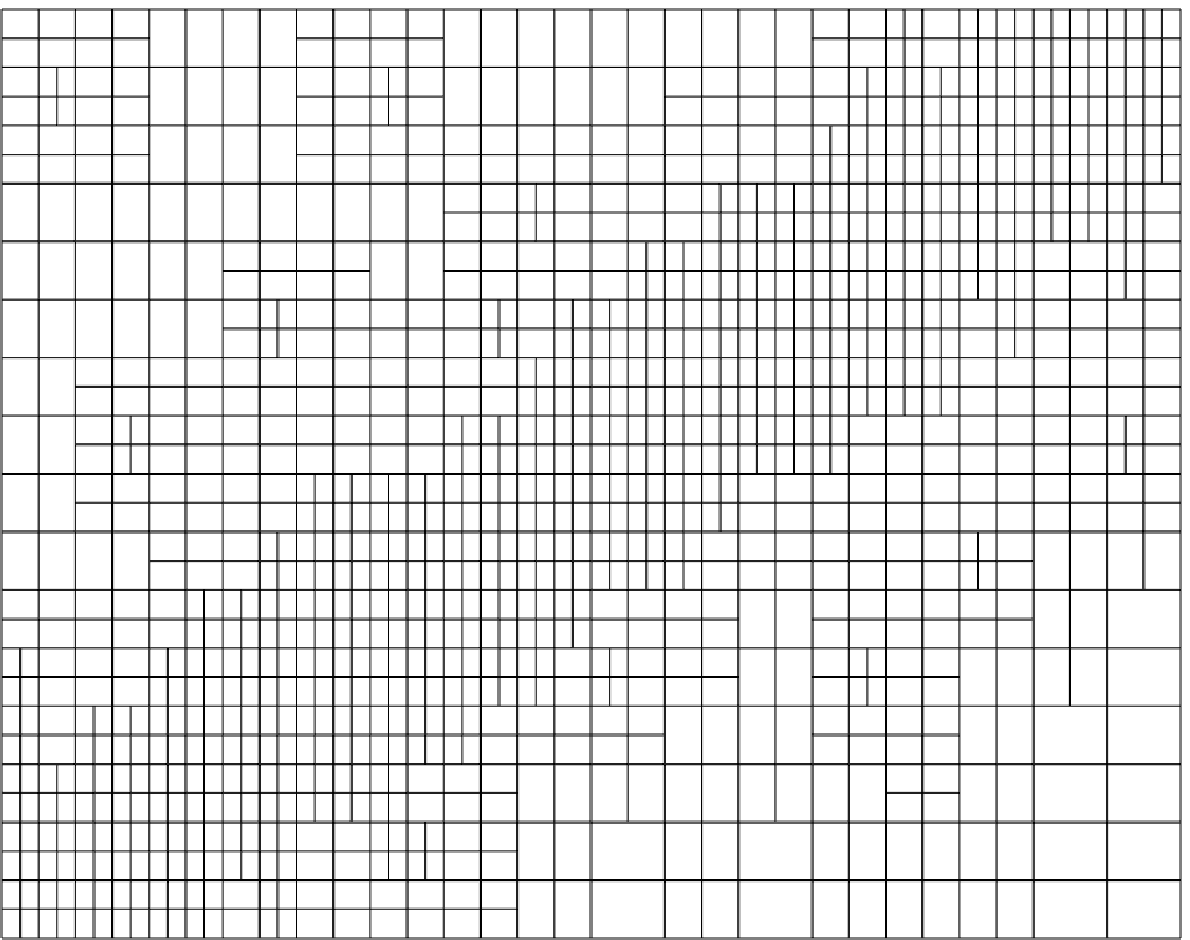}\hspace{.2\textwidth}%
\includegraphics[width=.3\textwidth]{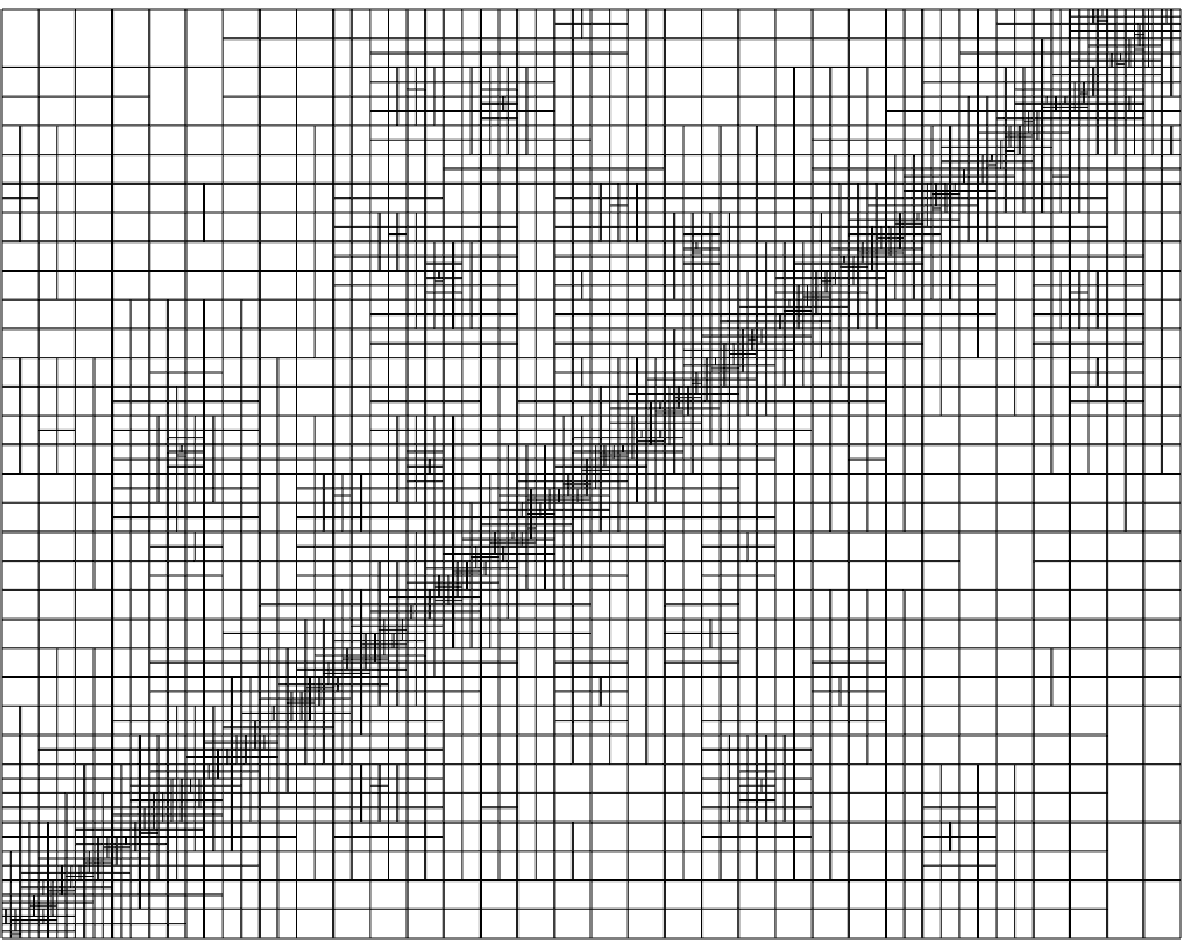}\\
\includegraphics[height=.225\textheight]{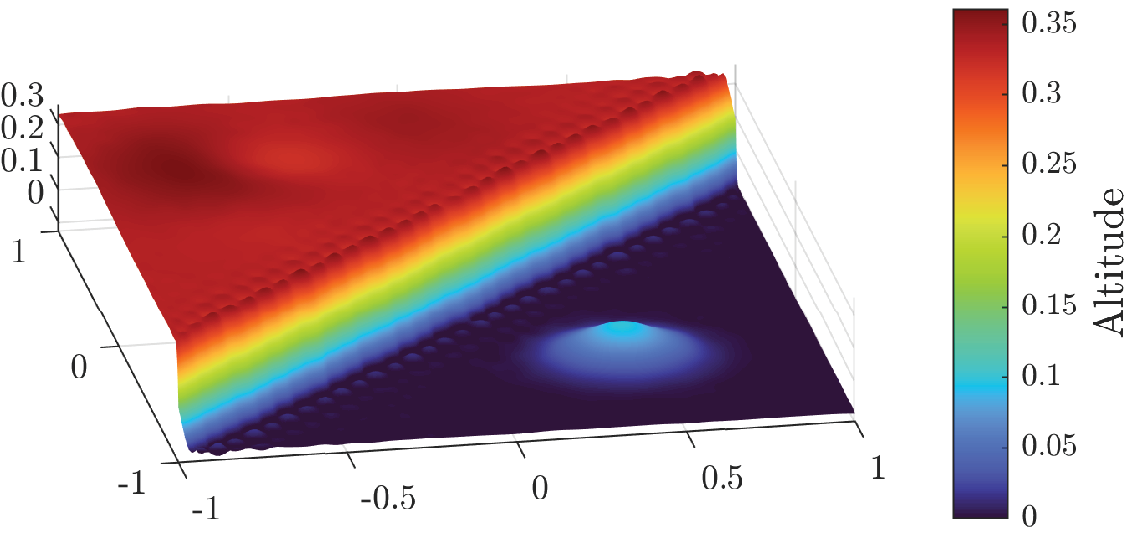}%
\includegraphics[height=.225\textheight]{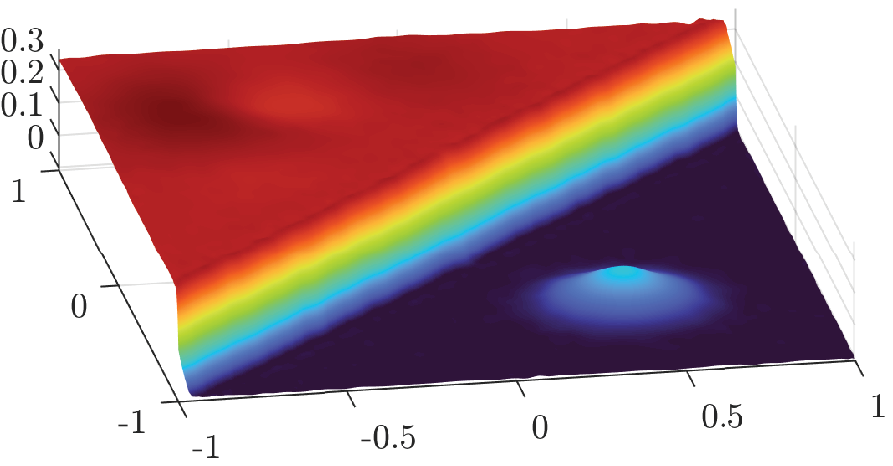}
\caption{Sharp surface. Top: T-mesh. Left: LS-T, right MTA. Bottom: visualisation of the fitted surface.}
\label{fig: sharp surface}
\end{figure}

\paragraph{T-meshes and fitted surface}
The T-mesh for the MTA strategy shows that more refinement is performed than for the LS approach, particularly in the dam domain. The reason is that the LS method minimizes the $\rmse$, which is also used pointwise as an error indicator for refinement. The MBA approach yields a greater $\rmse$ for the same refinement level and, hence, implies more refinement than the LS method. The consequence is a fitted surface $\Ssharp$ which oscillates much less than the LS approximation in the vicinity of the dam (\autoref{fig: sharp surface}, bottom left). This is the clear drawback of the global LS approach regarding the local explicit MTA. Please note that the point clouds were parametrized with the uniform method for the sake of simplicity. Other methods could have improved the results slightly. However, such investigations are beyond the scope of this contribution; see \cite{CampenZorin2017} for the use of similarity maps within the context of T-splines approximation. They do not affect the smoothness property of the MTA approach and/or its capacity to provide well-behaved approximation in the presence of data gaps and outliers, see \autoref{sec: outliers}.
\paragraph{Performance indicator}
Nine iterations are sufficient here with the MTA strategy to reach the same performance indicators as the LS-T with eight iterations. 
The computation time is slightly increased with ten iterations, but remains acceptable (approximately 1000 s for a non-optimised algorithm using Matlab on the computer mentioned previously). 
The $\rmsenoisek$  or $\rmsemathk$ are 0.0046 and 0.0036, respectively. 
The number of points outside the tolerance is more than twice smaller after ten iterations than after eight iterations for LS-T. We notice that for a comparable number of control points, the MTA gives a slightly lower RMSE than LS-T.
Consequently, sharp geometries can also be approximated with MTA, as the explicit local strategy acts indirectly to improve the parametrization locally.

\subsubsection{Outliers and data gaps}\label{sec: outliers}
\begin{figure}[ht]
\centering 
\includegraphics[height=.225\textheight]{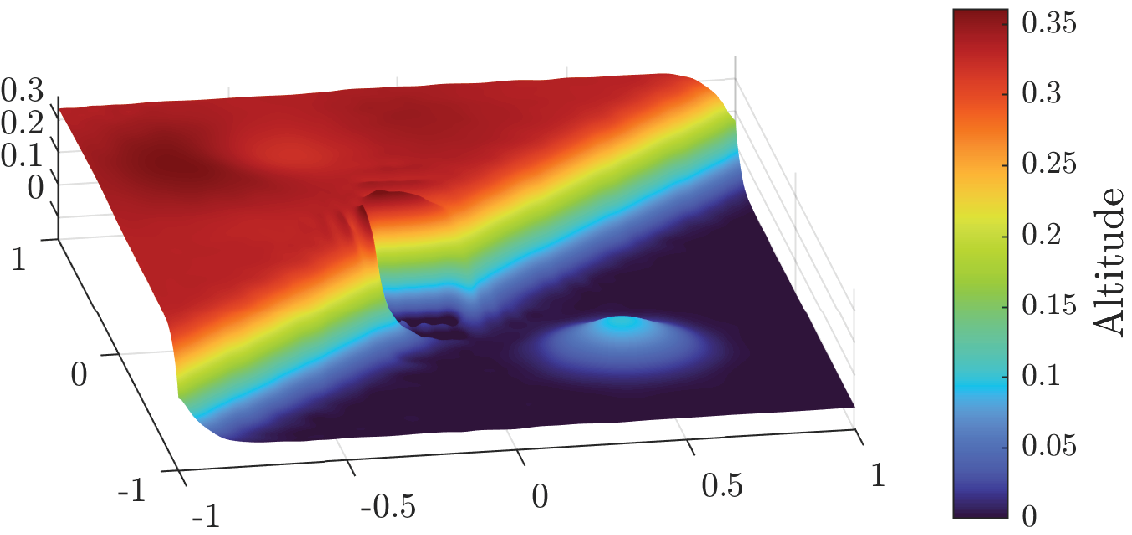}%
\includegraphics[height=.225\textheight]{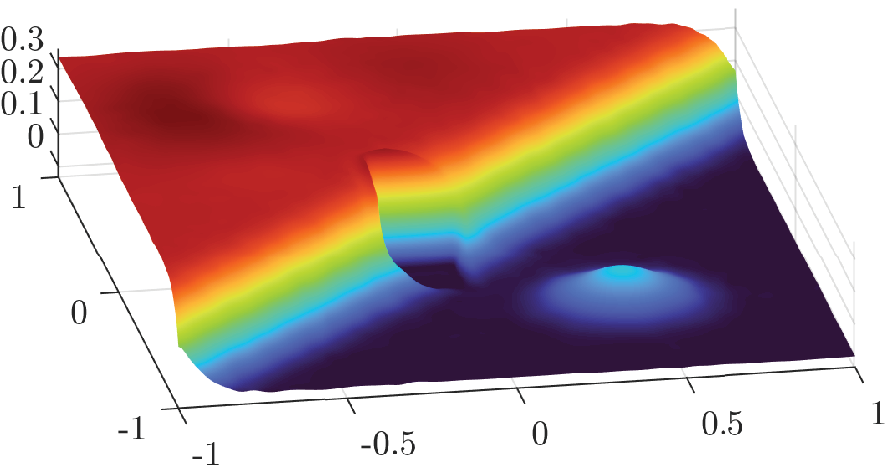}
\caption{Smooth surface with a data gap, approximated with LS-T (left) and MTA (right). 
}
 \label{fig: surface with gaps}
\end{figure}

\paragraph{Data gaps}
The results of the previous section are confirmed in the presence of a data gap. Both the T-meshes and the fitted surfaces presented in \autoref{fig: surface with gaps} show the high support provided by the MTA strategy regarding the LS-T near the gap. Correspondingly, oscillations and ripples occur in the LS-T fitted surfaces, which propagate far away from the gap itself. It yields a high maximal error of \num{.2}, i.e. four times higher than for the MTA. This is not the case for the MTA approximated surface: combined with a two times smaller $\rmse$ regarding the LS-T, a lower number of points outside the tolerance for the same computational time and a higher degree of details, the advantages of MTA are substantial. When the same number of controls points are estimated with the two strategies, i.e. disregarding the increase of computation time with LS-T, the $\rmse$ with LS-T is 0.001 higher than with MTA, and the maximum error 0.05. 
The mathematical proof why MTA performs better for point clouds with data gaps remains the topic of a subsequent contribution; the authors conjecture that this is due to the closeness of the explicit MTA strategy with an L1-norm minimization. 

\paragraph{Outliers}
We kept the same $\thr$ to fit $\Soutliers$ for the sake of comparison with the results obtained with $\Ssmooth$. An increase of the threshold or a fairness function could meet the challenge of overfitting and would smooth the fitted point clouds. The MTA and LS-T perform similarly, which we expect from $\Ssmooth$. The number of points outside tolerance is nearly nine times higher, the $\rmse$ three times and $\ncpk$ twice compared to $\Ssmooth$. These results are justified by the fact that the algorithm tries to fit the noise, although this effect should be limited by the adaptive refinement.

\begin{table}[ht]
\centering
\begin{tabular}{clccrccr}
\rowcolor{headrow}\cellcolor{white} && $\rmsenoisek$ & $\maxerr[k]$ & $n_{\text{out},k}$ & $\rmsemathk$ & $\ncpk$ & $\ct_k$\\
&MTA & \num{0.0032} & \num{0.0173} & \num{128} & \num{0.0014} & \num{2879} & \num{752} \\
\cellcolor{white} &LS-T  & \num{0.0033} & \num{0.0169} & \num{243} & \num{0.0016} & \num{1265} & \num{771} \\
\multirow{-3}{*}{ smooth}
&NURBS & \num{0.0051} & \num{0.0238} & \num{33769} & \num{0.0019} & \num{2345} & \num{1250} \\
\cellcolor{white} &LS-T & \num{0.0030} & \num{0.0161} & \num{132} & \num{0.0014} & \num{2972} & \num{2508} \\
\hline 
\cellcolor{white} &MTA & \num{0.0046} & \num{0.0507} & \num{1058} & \num{0.0036} & \num{3563} & \num{1150} \\
\cellcolor{white} &MTA  & \num{0.0061} & \num{0.0615} & \num{1940} & \num{0.0054} & \num{2471} & \num{852} \\
\multirow{-3}{*}{sharp}
\cellcolor{white} &LS-T & \num{0.0065} & \num{0.0674} & \num{2532} & \num{0.0058} & \num{1480} & \num{1021} \\
\cellcolor{white} &LS-T & \num{0.0063} & \num{0.0623} & \num{2089} & \num{0.0057} & \num{2392} & \num{2670} \\
\hline 
&MTA & \num{0.0046} & \num{0.0511} & \num{1037} & \num{0.0036} & \num{3502} & \num{1495} \\

\cellcolor{white} &LS-T & \num{0.0083} & \num{0.2104} & \num{1162} & \num{0.0077} & \num{1197} & \num{1508} \\
\multirow{-2}{*}{data gaps}
\cellcolor{white} &LS-T & \num{0.0059} & \num{0.1051} & \num{1089} & \num{0.0045} & \num{3418} & \num{4025} \\
\hline 

\cellcolor{white}&MTA & \num{0.0121} & \num{0.4605} & \num{1018} & \num{0.0117} & \num{4534} & \num{1456} \\
\multirow{-2}{*}{outliers}
\cellcolor{white} &LS-T & \num{0.0126} & \num{0.4715} & \num{1238} & \num{0.0122} & \num{1674} & \num{1624} \\

\cellcolor{white} &LS-T & \num{0.0122} & \num{0.4621} & \num{1031} & \num{0.0119} & \num{4489} & \num{5008} \\
\end{tabular}
\caption{Results of the approximation of the LS-T, MTA and NURBS strategies. The number of iterations is fixed to \num8, \num8 and \num{10} for the LS-T, NURBS and MTA, respectively. The performance indicators are those specified in \autoref{sec: performance indicator}.}
\label{tab: approximation results}
\end{table}

\subsection{Conclusion}
We simulated a high number of point clouds corresponding to a real case scenario. The following conclusions could be drawn:
\begin{itemize}
\item The MTA and LS-T strategies perform similarly for a smooth geometry, provided that the iteration number is higher for the MTA compared to the LS-T. We further mention that the number of iteration steps is left to the user’s convenience; this value should be a balance between the computation time and saturation of the $\rmse$ to avoid noise fitting for a predefined threshold. At that point, there is no way to give a rule of thumb. 
\item The surface approximation using NURBS and global LS was shown to perform suboptimally. Here, the local refinement strongly avoids noise fitting in domains where no refinement is needed after a few iteration steps. 
\item The MTA outperforms the LS-T in the presence of outliers, data gaps or for sharp geometries. The approximation of the residuals is known to act as a smoothing: it avoids polluting effects, i.e. propagating ripples and oscillations originating from the missing points along the whole surface.
\item The MTA is an explicit approach after a few iterations with a global LS. This results in a strong decrease of the computation time compared to LS-T for a given iteration step as no matrix inversion or product is needed. A higher degree of refinement and the computation of more iterations can easily be done to fit the surface optimally. 
\end{itemize}

\section{Real data analysis}\label{sec: real data}
In this section, we propose fitting real point clouds of TLS observations. We approximate two kinds of surfaces for which no LS-T solution could be computed due to the unfavourable data density of the point cloud: a 3D print of the mathematical surface defined in \autoref{sec: point clouds} and a sand-dune. We aim to highlight the potential of T-splines for new applications for land remote sensing to identify deformation within a GIS context, i.e. beyond the CAD framework. An implementation of the MTA method in CAD or GIS software is easy and could done at a low cost.
\subsection{Real dataset 1: 3D printed surface}
The mathematical surface defined in \autoref{fig: smooth surface} was printed in 3D and scanned under different scanning configurations by varying the tilting of the surface and the position of the TLS. We selected two point clouds of the MTA surface fitting for further analysis. Unusually, the reference is known exactly here -- up to the artefacts introduced by the 3D printer. This is rarely the case in real applications. This allows for a comparison with the true solution, unknown in real applications.
\subsubsection{3D printed surface: generation and scanning}
\begin{figure}[ht]
\centering
\includegraphics{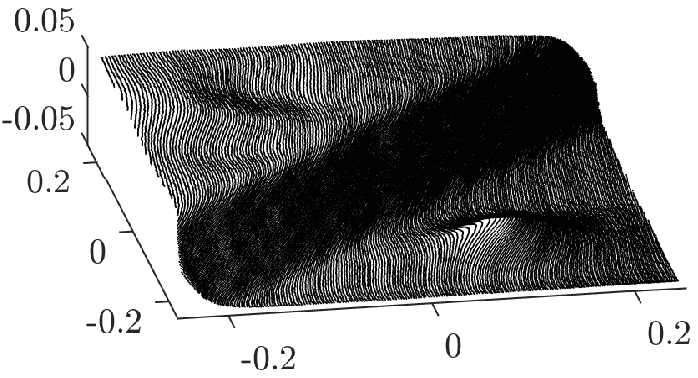}%
\includegraphics{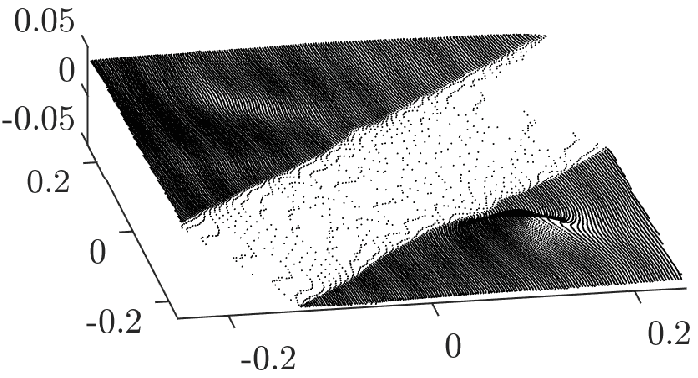}\qquad\\
\includegraphics{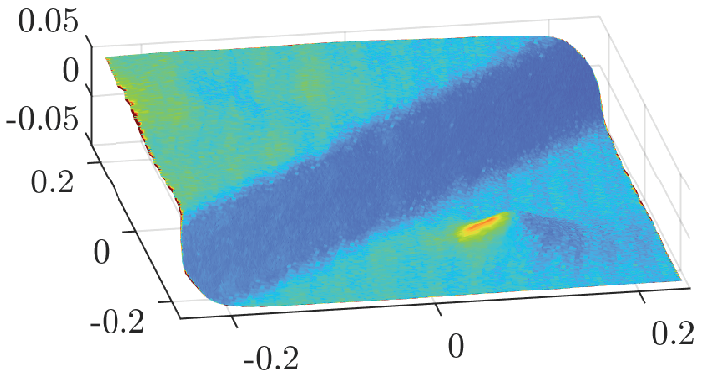}%
\includegraphics{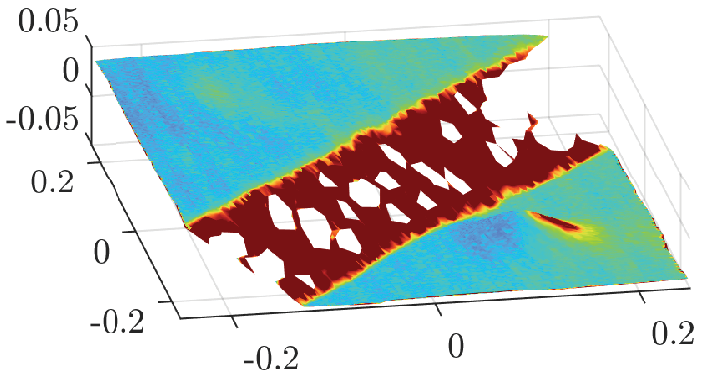}~
\includegraphics{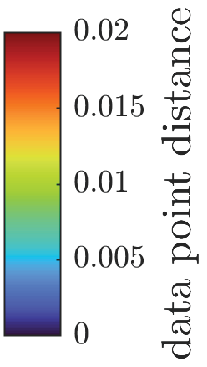}\\
\includegraphics[height=4cm]{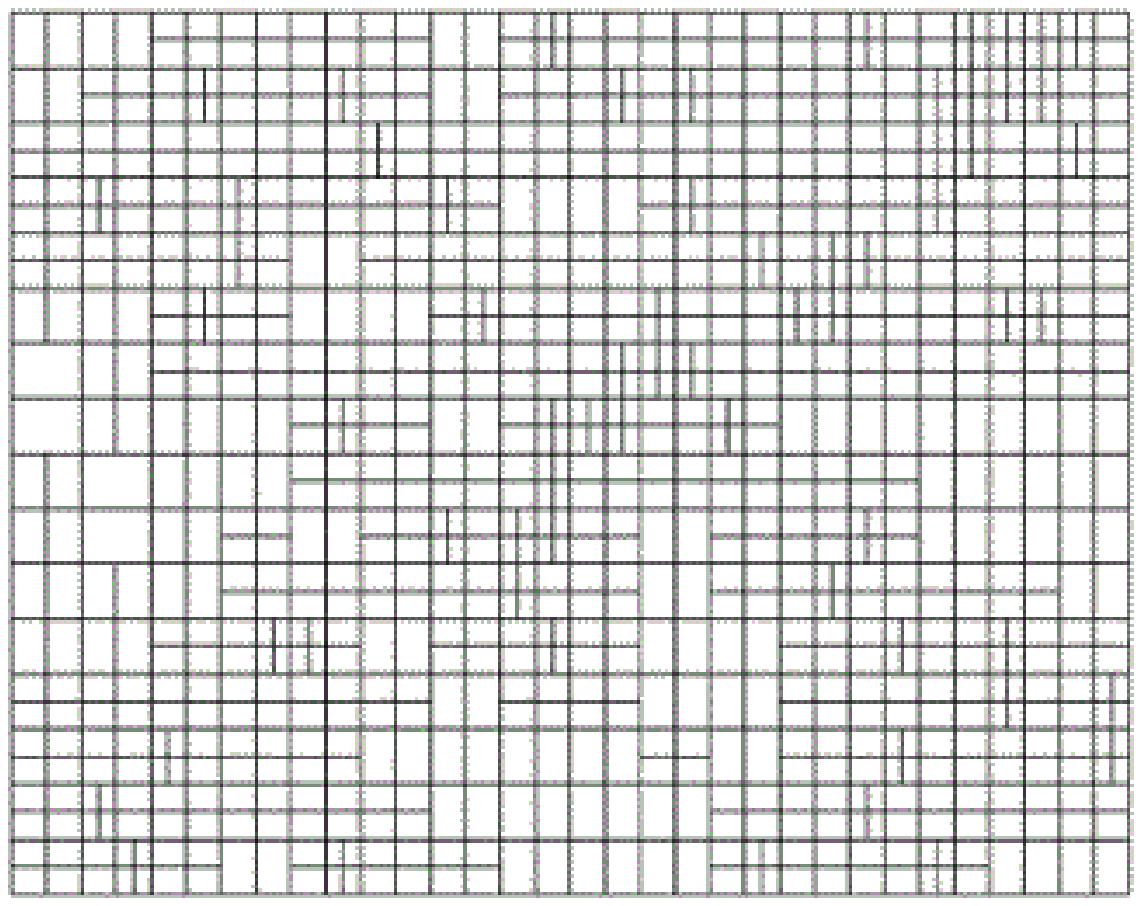}
\includegraphics[height=4cm]{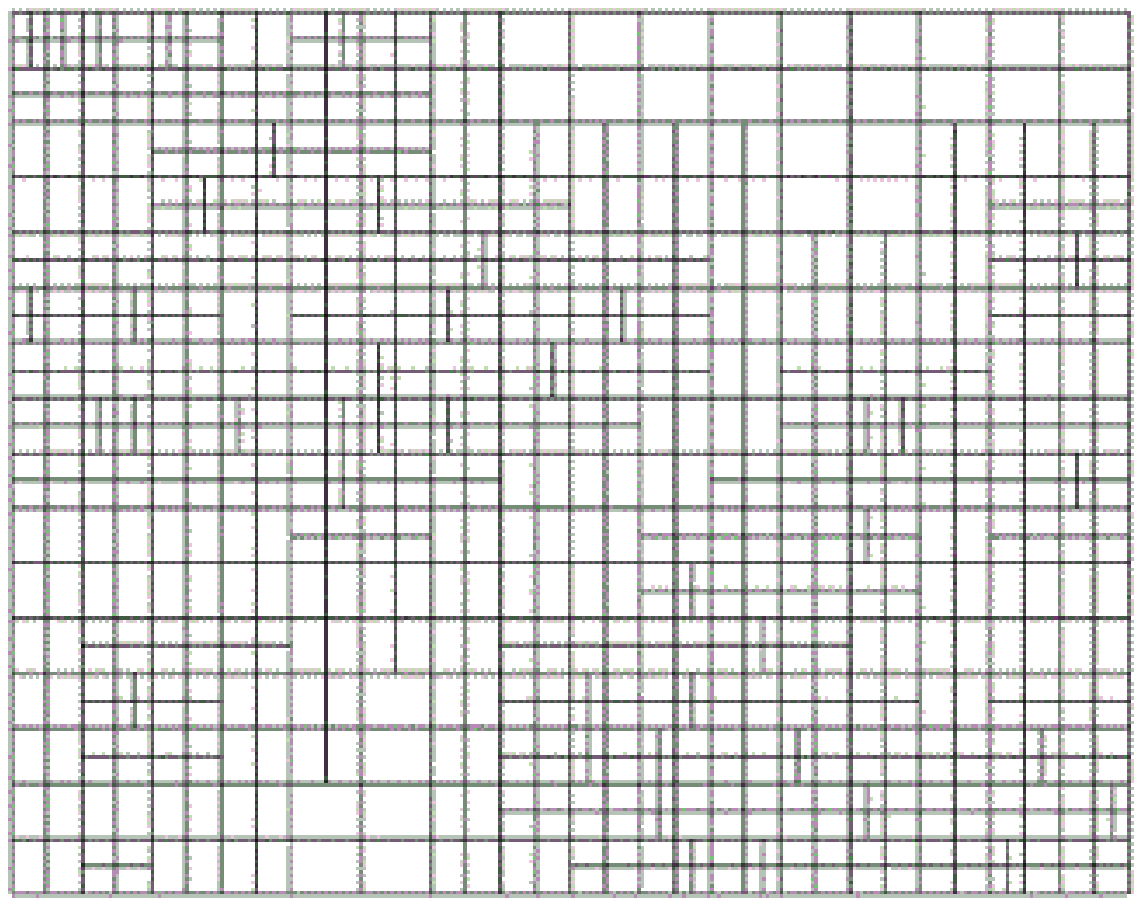}%
 \caption{Left: data from the first scanning configuration. Right: data from the second, unfavourable configuration. The black dots (top) represent the data points, the color scale (middle) represents the distance between data points. Bottom: the corresponding T-meshes}
 \label{fig: scanning configuration}
\end{figure}

The 3D printing of a simulated point cloud corresponding to \autoref{sec: simulations} is described in \autoref{apx: 3d printing}. 
The measurements of the surface took place in the measuring laboratory of the Geodetic Institute, Leibniz University Hannover, with a Zoller+Fröhlich Imager 5016 phase scanner (see \autoref{fig: scanning configuration}, left). The position of the laser scanner was varied so that shadowing occurred; the latter is necessary as it allows testing the MTA fitting methodology in case of outliers and data gaps, as mentioned in \autoref{sec: simulations}. The outliers and data gaps occur along the dam and for the mountain-like Gaussian bell. Please refer to \cite{KermarrecSchildHartmann2021} for more details on the experiment.
In this contribution, we make use of two point clouds scanned under a suboptimal position of the TLS at a distance of approximately \SI 3\m. There are approximately \num{57 680} points for the first point cloud and \num{97490} for the second one. The approximation remains manageable in a standard computer using Matlab. We did not measure the tilting exactly or the exact distance for the scanning configuration: we wanted to test the feasibility of using the MTA in a real scenario and not study the shadowing regarding the scanning configuration. 
The two point clouds under consideration are visualized in \autoref{fig: scanning configuration} (middle and right): data gaps are clearly visible. The threshold chosen for the approximation was fixed at \SI{0.01}{\mm}, following the manufacturer's specifications for the expected noise level in the $z$-direction. The design matrix became rapidly singular due to the missing observations. Matlab solvers reach their limit for this unfavourable data density. The MTA strategy here can easily provide a surface fitting from the third iteration using adaptive local refinement with T-splines. Unfortunately, the Matlab software ran out of memory for ten iterations due to the high number of points to process; we, therefore, restricted ourselves to a maximum of eight iterations. 

\subsubsection{Results}
The results of the MTA fitting are summarized in \autoref{tab: approximation results}. 
The same performance indicators are used here as those for the simulations.

\paragraph{T-meshes}
The T-meshes highlight the domains in which no points were available for the surface fitting, i.e.\ the dam for the first configuration and the bell for the second one: no refinement can be performed in that domain due to the data density. As expected, small artefacts coming from the 3D printing are fitted additionally; the result is additional control points to fit these ``outliers''. Thus, we do not expect the T-mesh to be more refined in the dam or Gaussian bell domain.
\paragraph{Performance indicator}
\begin{table}[ht]
\centering
\begin{tabular}{lcccccc}
\rowcolor{headrow} Real point clouds & $\rmsenoisek$ & $\maxerr[k]$ & $n_{\text{out},k}$ & $\rmsemathk$ & $\ncpk$ & $\ct_k$\\
First  configuration & \num{0.0030} & \num{0.0131} & \num{45} & \num{0.0007} & \num{1086} & \num{211} \\
Second configuration & \num{0.0029}  & \num{0.0133} &  \num{86} & \num{0.0003} & \num{1283} & \num{376}
\end{tabular}
\caption{Real case scenario using the MTA strategy with $\thr=\SI{0.01}{\mm}$ and a maximum of eight iterations. The results are given in \si{\mm}.}
\label{tab: real case scenario}
\end{table}

\autoref{tab: real case scenario} highlights the high trustworthiness of the approximation regarding the original mathematical surface: $\rmsemathk$ is below the \si{\mm} level. The number of points outside the tolerance was low for the first configuration considering the threshold chosen and a maximum number of iterations. This configuration was slightly more favourable than the second one, for which a strong shadowing occurs. Nevertheless, the  $\rmsenoisek$ reaches \SI 3 \mm, which is comparable to the noise level expected for objects scanned at short distances. This example shows that the MTA can be used to approximate challenging point clouds with varying data density where a global LS strategy would fail.
\subsection{Sand-dune}
\subsubsection{Description of the dataset}
\begin{figure}[ht]
\centering 
\includegraphics{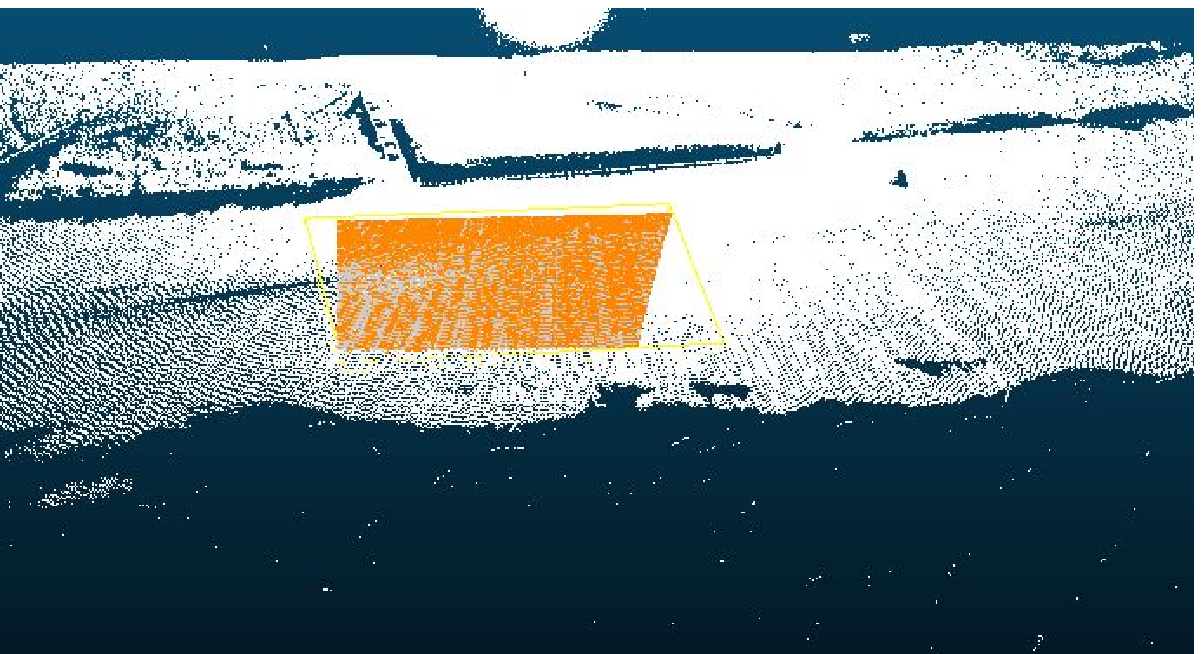}\\[1ex]
\includegraphics[width=.45\textwidth]{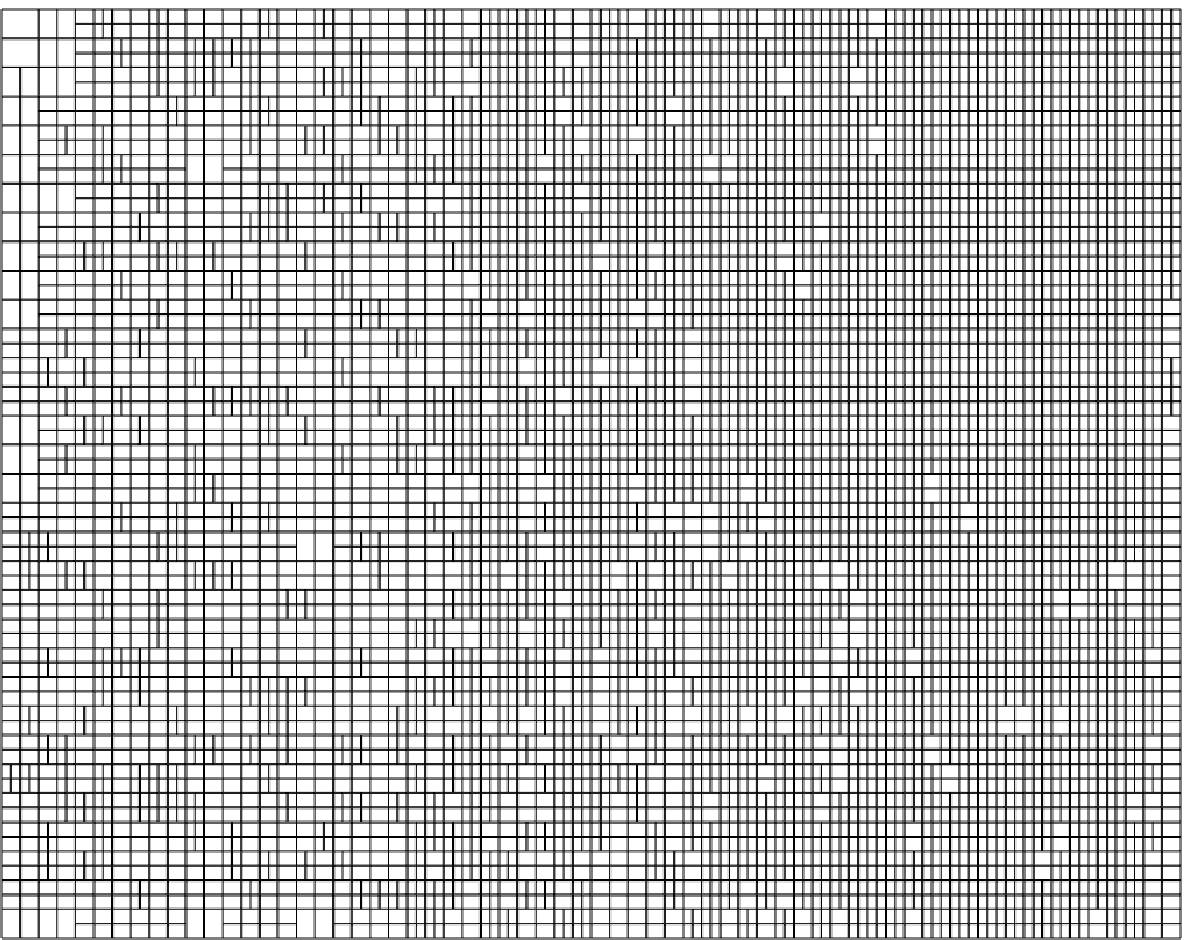}\hspace{.05\textwidth}%
\includegraphics[width=.45\textwidth]{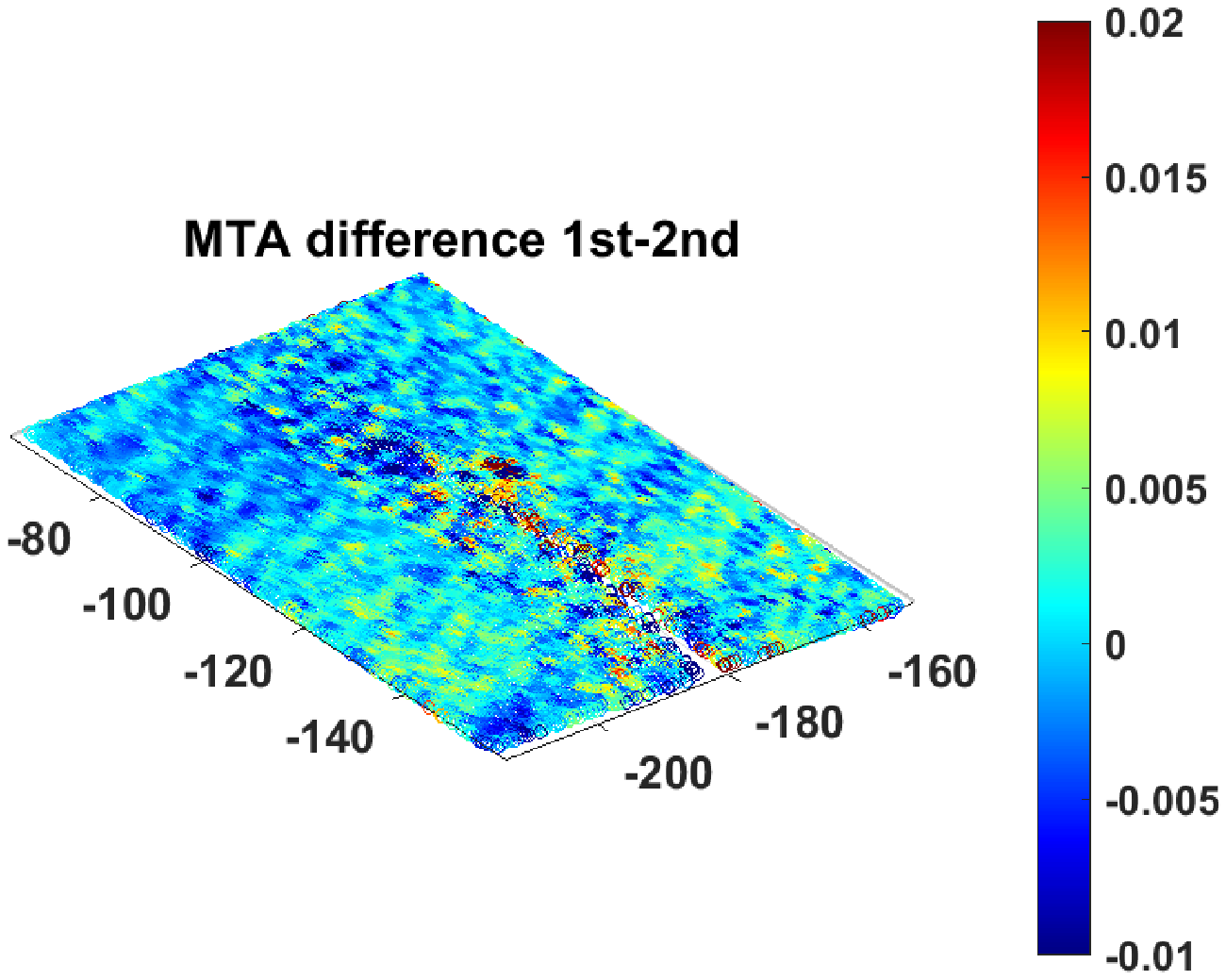}
 \caption{Top: localisation of the point cloud of the dune selected. Bottom: the T-mesh (left) and surface difference (right) after ten iterations using the MTA strategy for the first point cloud with  $\thr=\SI{0,01}{\m}$.
 }
 \label{fig: mta at the beach}
\end{figure}

In order to further demonstrate the potential of the MTA method, we made use of a dataset%
\footnote{available at \url{https://data.4tu.nl/articles/dataset/CoastScan\_Data\_of\_daily\_scans\_at\_low\_tide\_Kijkduin\_January\_2017/12692660}} 
that was acquired for the CoastScan project \cite{VosKuschnerus2020} along the coast of Kijkduin in the Netherlands. A Riegl VZ2000 laser scanner was mounted on the roof of a hotel and programmed to perform a scan of the nearby dune and beach area every hour. In this contribution, we made use of two datasets scanned on consecutive days in January 2017. The corresponding .asc files are \texttt{170101\_230044} and \texttt{170102\_010046} and can be freely downloaded. The files were corrected for a tilt in the scanner if the tilt exceeded \num{.01} degree of the median inclination angle recorded. More information about the preprocessing of the point clouds can be found in the dedicated publication. Using the free software Cloud Compare, we cut a small nearly rectangular surface of the beach scanned, illustrated in \autoref{fig: mta at the beach} (top). There were approximately \num{21 000} points for each epoch. 

\subsubsection{Fitting with MTA and visual deformation analysis}
We intentionally present the results of the performance indicators for the first epoch only, as the second epoch gave similar results, see \autoref{tab: real case scenario} and \autoref{fig: mta at the beach}. The T-mesh after ten iterations highlights the level of refinement performed and the potential of the MTA to fit point clouds with a data gap (see \autoref{fig: mta at the beach}, top). The $\rmsenoisek[10]$ after the tenth iteration was \SI{0.01}{\m} with \num{7219} control points computed and a relatively high number of points outside the tolerance (\num{2848}). We know from the simulations that this is an indication that the threshold chosen may be too low regarding the noise level. This finding is supported by the fact that the stagnation level for the $\rmse$ was not reached after ten iterations (the corresponding results are not shown here for the sake of brevity but can be made available on demand). Consequently, we increased the threshold to \SI{0.02}{\m}. After ten iterations and only \SI{230}{\s}, the $\rmsenoisek[10]$ was similar to the one obtained with a threshold of \SI{0.01}{\m}. The number of points outside the tolerance decreased from \num{2848} to \num{477} for the same $\maxerr[,\text{noise},10]$. This higher threshold seems to be more adequate to fit the point cloud under consideration with a high trustworthiness by simultaneously decreasing the risk of overfitting.

\begin{table}[ht]
\centering
\begin{tabular}{lccrccr}
\rowcolor{headrow} Threshold & $\rmsenoisek$ & $\maxerr[k]$ & $n_{\text{out},k}$ & $\rmsemathk$ & $\ncpk$ \\		
$\thr=\num{.01}$ & \num{0.0256} & \num{0.9410} & \num{2848} & \num{7219} & \num{2170} \\
$\thr=\num{.02}$ & \num{0.0259} & \num{0.9374} & \num{477} & \num{1043} & \num{234} \\
\end{tabular}
 \caption{Performance indicator for the first epoch \texttt{170101\_230044} and two thresholds.}
 \label{tab: performance indicator}
\end{table}

In order to highlight the potential of the adaptive approximation for performing a deformation analysis with T-splines surfaces, we estimated the parametric surfaces from the two epochs at the same horizontal coordinates.
This is made possible by the mathematical expression of the surface and is one of the major advantages of surface fitting for deformation analysis, see \cite{KermarrecSchildHartmann2021}.
We took the difference of the $z$-component between the gridded point clouds. The result is depicted in \autoref{fig: mta at the beach} (bottom right). The difference reaches a maximum at the cm level and could provide some indications about the movement of the sand-dune on a daily basis, provided that this level does not exceed the unknown registration error. The differences near the shadowed area are slightly stronger (from \num{-1} up to \SI 2 \cm).
An exact analysis of the causes remains beyond the scope of the present contribution, which aimed to show the potential of surface fitting with MTA for applications in land remote sensing in a GIS context only.
\section{Conclusion}\label{sec: conclusion}
The fitting of huge scattered and noisy point clouds of observations with variable data density from contactless sensors with parametric surfaces is a challenging task. Mathematical approximations should 
\begin{inlinelist}
 \item avoid the fitting of the noise 
 \item be computationally manageable 
 \item provide trustworthy and interpretable results for GIS applications, such as deformation analysis. 
\end{inlinelist}
Unfortunately, the NURBS surface does not allow for local refinement due to the tensor product formulation. The associated drawbacks can be addressed by means of adaptive local refinement strategies. We focused on T-splines in this contribution. This type of spline enables local knot refinement to avoid superfluous control points by allowing T-junctions in the control mesh, i.e. a rough approximation of the surface itself. The T-splines are fully compatible with NURBS surfaces. The latter are a special case of T-spline surface and have been widely accepted in the design communities and applied to different research areas. Their potential to fit observations from TLS for deformation analysis within a GIS context have, to the best of the authors’ knowledge, not been assessed up to now.
The iterative fitting of point clouds is often performed with an LS adjustment, which necessitates matrix product and inversion. When a high number of control points have to be estimated or the points on the surface are not close enough to the grid structure, the design matrix rapidly becomes singular, forbidding further refinement. We have introduced a new method in this contribution. Our fitting procedure combines the AST-splines refinement with a multilevel approach based on approximating the residuals of the surface fitting. This strategy, called the MTA, was shown to be highly trustworthy in an MC framework, particularly in the case of data gaps. Oscillations and ripples were avoided for sharp geometries and in the presence of outliers. The slight increase of parameters to estimate was counterbalanced by the computational efficiency of the method, which is very simple to implement as explicit, stable and fast. These advantages pave the way for its wider use in the geodetic context for land sensing and GIS applications, such as deformation analysis. We have applied the newly developed method to real point clouds for which no LS approximation could be performed due to shadowing effects and irregular and non-rectangular point clouds. We printed the mathematical surface used in the simulations in 3D and scanned it with a TLS under an unfavourable scanning configuration. We could show that the $\rmse$ for high shadowed point clouds regarding the scanned surface was under the mm level concerning the mathematical surface. The computation time was manageable on a standard computer. This result highlighted the  trustworthiness of the MTA regarding the LS strategies. The latter may lead to surfaces which are not well-behaved due to numerical inaccuracies as the level of refinement increases. We confirmed these results using a point cloud from a sand-dune scanned in the Netherlands. Further investigations will be performed by printing a more challenging surface to fit in 3D. Local threshold strategies should be implemented with the aim of improving the MTA fitting further.

\section{Acknowledgments}
The authors warmly thank Niklas Schild for having implemented the MTA successfully and performed the measurements. The Nederlandse Organisatie voor Wetenschappelijk Onderzoek is thanked for having founded the project which led to the freely available point clouds of the sand-dune. This study is supported by the Deutsche Forschungsgemeinschaft under the project KE2453/2-1.

\appendix

\section{3D printing and measurements}\label{apx: 3d printing}
The 3D printing of a simulated point cloud starts with a digital volume model, which is converted into a STL 
file compatible with the 3D printer. In this format, the surface of the object is described by a network of 
triangular areas and the object divided into layers. A common and widespread 3D printing method is fused 
layer modelling; We used a polyactid (PLA-HAT) which is melted by a heating element on the printhead and 
applied to the location desired using the movement possible in three spatial axes. The corpus printed has 
a size of $50{\times}\SI{50}{\cm}$, limited by the specificity of the 3D printer used: we aimed to print 
the surface without junctions to avoid unnecessary artefacts. Such effects are unavoidable even if the 
surface is printed as a whole, as they come from a post-printing polishing or droplets lost by the 3D printer.
They can be considered as ``outliers'' regarding the reference surface. In a last step, the 3D surface obtained
was lacking in grey tones so that no strong reflectivity occurred during the scanning process. 
The measurements of the surface took place in the measuring laboratory of the Geodetic Institute, Leibniz
University Hannover, with a Zoller + Fröhlich Imager 5016 phase scanner (Zoller \& Fröhlich GmbH, Wangen, 
Allgäu). This instrument is well-designed for short-range indoor measurements. We varied the position of the
laser scanner so that shadowing occurred; the 3D-printed planes were cut out of the whole point cloud using 
the free software CloudCompare \cite{CloudCompare} to process the 3D point clouds. The edges of the plane 
were treated carefully to avoid mixed pixels. We used the TLS settings ``quality high'' and 
``resolution high'', which resulted in approximately \num{57 680} points for the first point cloud
 and \num{97 490} for the second one.


\end{document}